\documentclass[12pt,a4paper,twoside,final,notitlepage, leqno]{article}
\usepackage[english]{babel}
\usepackage[T1]{fontenc}  
\usepackage{graphicx}
\usepackage{amsmath,amsthm,epsfig,amsfonts,bbm}  

\newcommand{\pequationdeb}{$$ \left\{ \begin{minipage}[c]{130mm}}
\newcommand{\pequationfin}{\end{minipage}
                           \right. $$}

\def \smb {{\scriptstyle \bullet }}
\newcommand{\monitem}{ \smallskip \noindent $\bullet$ \quad  } 
\newcommand{\moneq}{\vspace*{-6pt} \begin{equation} \displaystyle } 
\newcommand{\moneqstar}{\vspace*{-6pt} \begin{equation*} \displaystyle } 
\newcommand{\monendstar}{\vspace*{-6pt} \end{equation*}   }
\newcommand{\monend}{\vspace*{-6pt} \end{equation}   }
\newcommand{\beq}     {\begin{equation}}
\newcommand{\enq}     {\end{equation}}
\newcommand{\be}    {\begin{enumerate}}
\newcommand{\ee}    {\end{enumerate}}

\newcommand{\Bb}

\def\R{{\rm I}\! {\rm R}}

\def\section*#1{}
\def\resume{\if@twocolumn
\section*{R\'esum\'e}
\else \small
\quotation{\bf \it R\'esum\'e \rule[1mm]{1.5mm}{0.2mm}\vspace{0pt}}
\fi}
\def\endresume{\if@twocolumn\else\endquotation\fi}
\def\abstract{\if@twocolumn
\noindent\section*{{\bf Abstract}}
\else \small
\quotation{\noindent \bf {Abstract.} \rule[1mm]{1.5mm}{0.2mm}\vspace{0pt}}
\fi}
\def\endabstract{\if@twocolumn\else\endquotation\fi}

\usepackage{fancyhdr}
\renewcommand{\headrulewidth}{0pt}
\begin{document}

\fancypagestyle{plain}{ \fancyfoot{} \renewcommand{\footrulewidth}{0pt}}
\fancypagestyle{plain}{ \fancyhead{} \renewcommand{\headrulewidth}{0pt}} 

~
 
~

\bigskip \bigskip   \bigskip

\centerline {\bf \LARGE  Stable  lattice Boltzmann schemes}

\bigskip    \centerline  {\bf    \LARGE with a dual entropy approach}

 \bigskip  \centerline  {\bf     \LARGE for monodimensional nonlinear waves} 

\bigskip \bigskip \bigskip

\centerline { \large  Fran\c{c}ois Dubois~$^{ab}$}   

\bigskip  

\centerline { \it  \small    $^a$   Conservatoire National des Arts et M\'etiers, } 
\centerline { \it  \small   Department of Mathematics and LMSSC,  Paris, France.}  
\centerline { \it  \small  $^b$   Department of Mathematics, University  Paris-Sud,}  
\centerline { \it  \small  B\^at. 425, F-91405 Orsay Cedex, France.} 
\centerline { \it  \small francois.dubois@math.u-psud.fr}    \bigskip   
 
\bigskip


\centerline {  20  September 2012~\footnote{~Contribution published  in  
  {\it Computers And Mathematics With Applications},  volume~65, pages~142-159, 2013, 
doi: 10.1016/j.camwa.2012.09.005. Edition 04 March 2013. }}
 
\bigskip

\bigskip 
\noindent  {\bf Abstract. } \qquad 
We follow the mathematical framework proposed by Bouchut  \cite {Bo03}
and present  in this contribution a  dual entropy approach for determining equilibrium states of a 
lattice Boltzmann scheme. This method  is expressed in terms of  the
dual of the mathematical entropy relative to  the underlying conservation law.
It  appears as  a good mathematical framework for 
establishing a ``H-theorem'' for the system of equations with discrete velocities. 
The dual   entropy approach is used with D1Q3 lattice Boltzmann 
schemes for the Burgers equation. It conducts to the explicitation 
 of three different  equilibrium distributions of particles 
and induces naturally a nonlinear  stability condition. 
Satisfactory  numerical results for 
strong nonlinear shocks and rarefactions are presented. 
We prove also that   the dual entropy approach
can be  applied with a D1Q3 lattice Boltzmann scheme
for  systems of linear and  nonlinear acoustics and 
we present a numerical result  with strong nonlinear waves 
for nonlinear acoustics. 
We    establish also a negative result: 
with the present framework, 
the  dual entropy approach cannot be used for   the shallow water equations.                
 
 $ $ \\  [2mm]

\bigskip \bigskip  \newpage \noindent {\bf \large 1) \quad  Introduction}  

\fancyfoot[C]{\oldstylenums{\thepage}}
\fancyhead[OC]{\sc{Stable  lattice Boltzmann schemes with a dual entropy approach}}

\quad 
An hyperbolic partial
  differential equation like the Burgers equation 
\moneq  \label{abs-01}   
 \partial_t  u   \,+\,   \partial_x  \big(   F(u) \big) \,=\, 0  
 \,, \quad F(u) \, \equiv \,  {{ u^2}\over{2}} \, 
\monend   
exhibits shock waves  (see {\it e.g.} \cite{GR96}), 
{\it id est} discontinuities propagating with finite  velocity. 
 In order to select the physically relevant weak  solution, 
it is necessary to enforce the so-called entropy condition
\moneq  \label{abs-02}  
\partial_t \big( \eta(u) \big)  \,+\,  \partial_x  \big( \zeta(u) \big)  \,\leq \, 0   \,
\end{equation}      

\noindent 
as suggested by Godunov \cite{Go61} and 
Friedrichs and Lax \cite{FL71}.  
In the relation (\ref{abs-02}), 
$\, \eta(\smb) \,$ is a strictly convex function and $\, \zeta(\smb) \,$ the associated entropy
flux  (see  {\it e.g.} \cite{GR96},   \cite{DD05} or \cite{La06}). 
For the Burgers equation,  the quadratic entropy is usually considered   
\moneq  \label{quad-entrop}  
\eta(u)  \equiv  {{ u^2}\over{2}} \,, \qquad 
   \zeta(u)   \equiv   {{ u^3}\over{3}}  \, . 
\monend             
%
Remark that the entropy condition  (\ref{abs-02})  is just one of at
least three possible criteria for selecting the physically relevant weak solution. One may
also consider the vanishing viscosity limit, or the Lax entropy criterion (see 
{\it e.g.}  \cite{GR96}   or \cite{La06}). 
%

\monitem  
%
The computation of discrete shock waves with lattice Boltzmann approaches began with 
viscous Burgers approximations in the framework of lattice gaz automata
(see   Boghosian and  Levermore \cite{BL87}, Elton \cite{El91}, 
Elton {\it et al.} \cite{ELR93}).  
With the lattice Boltzmann methods  described {\it e.g.}
 by Lallemand and Luo  \cite{LL00}, first tentatives were proposed 
by  d'Humi\`eres \cite{DDH92}, Alexander 
 {\it et al.} \cite{ACCD92}, Qian and Zhou  \cite{QZ98}.  
%
The study of nonlinear scalar  equation with the help of the
lattice Boltzmann scheme has been emphasized by 
Buick  {\it at al.}  \cite{BBGG2k} 
for nonlinear acoustics. 
The approximation of the Burgers equation with a quantum variant of the 
method has been presented  by Yepez \cite{Ye02}. 
%
 A  D1Q2 entropic scheme for the one-dimensional viscous 
Burgers equation has been developed by Boghosian  {\it et al.}  \cite{BLY04} 
%
 and we refer to 
Duan and Liu  \cite{DL07} for 
the approximation of two-dimensional Burgers equation.
The extension for gas dynamics equations 
%
and in particular shock tubes problems 
is under study with {\it e.g.} 
%
the  works of  Philippi {\it et al.}, \cite{PHSS90},  
Brownlee {\it et al.} \cite  {BGL07},  
Nie, Shan and Chen \cite{NSC08},   
 Karlin and Asinari \cite{KA10},  
Chikatamarla and Karlin \cite{CK09}.

%

\monitem  
%
In this contribution, we experiment the ability of lattice Boltzmann schemes to approach
weak entropic solutions of hyperbolic equations.  
In such situations, the scheme  exhibits   some kind of vanishing viscosity limit. 
%
We start from the mathematical framework developed 
by Bouchut  \cite {Bo03} making the link between the finite volume method and kinetic
models in the framework of the BGK \cite{BGK} approximation. 
The key notion is the representation of the dual entropy
with the help of convex functions associated with 
the discrete velocities of the lattice.  
We call ``dual entropy approach'' the set of associated constraints 
for the equilibrium distribution. 
%
In section~2, we recall  this framework with emphasis on  the one-dimensional 
case and prove a continuous version of the ``H-theorem''. 
In section~3 we derive three
equilibria for a D1Q3 kinetic distribution associated with the lattice Boltzmann 
method applied to the Burgers equation. 
In section~4, we precise our  numerical D1Q3 scheme and make a simple link 
with the finite volume approach. 
We present   numerical experiments with nonlinear  Burgers waves in section~5.
%
In   section~6, we study  the ability of the dual entropy approach to determine
D1Q3 equilibria for systems of  linear  and nonlinear acoustics.
We study  the system of shallow water equations in Section~7. 
%


\bigskip \bigskip  \noindent {\bf \large 2) \quad  Kinetic representation of the dual
  entropy} 

\quad  
The  Legendre-Fenchel-Moreau duality is a classic notion defined 
when     we consider  a convex function $\, \eta(\smb)   \, $ 
of several  variables. We can apply the   duality transform 
that suggests that   convex function  $\, \eta(\smb)   \, $   
is parametrized    by the slopes of the  tangent planes. In other terms, 
\moneq 
 \eta^* (\varphi) \, = \, \sup_{W}  \Big(  
\varphi \, \smb \, W \,\,-\,\,  \eta(W)   \Big)  \, . 
\label{dual}  \monend 
The upper bound in the right hand side of relation (\ref{dual}) 
is   obtained   (when it is not on the boundary of the domain of
variation  of the state $W$) 
by solving the equation of unknown~$ \, W \, $: 
\moneq \label{deriv-entrop} 
 \eta' (W) = \varphi \, . 
 \monend
A first example is simply $ \,  \eta(w) \equiv   {\rm e}^w   \, $ 
at one space dimension. Then  
$  {\rm e}^w = \varphi , \, $ $ \, \eta^* (\varphi) =  
 \varphi \, \log \varphi  -   \varphi   \,  $  
and we recover in this way the fundamental tool to define the so-called
``Shannon entropy'' \cite {Sh48}. 

\monitem 
We can derive the dual function~:  
if $   \, \, {\rm d} \eta(W) \equiv \varphi \,\smb \, {\rm d}W  \,  $
then 
%
\moneq  \label{dual-etat}
     {\rm d} \eta^*(\varphi)   \,=\,  {\rm d} \varphi \,\smb \,  W   
\monend   
%
%
and the ``physical state'' $W$  is the Jacobian of the  dual entropy. 
In an analogous way, we can introduce (see {\it e.g.} 
 \cite{GR96},   \cite{DD05} or \cite{La06})  in the context of hyperbolic conservation laws 
\moneq  \label{hcl}    
 \partial_t  W   \,+\,   \partial_x  \big(   F(W) \big) \,=\, 0   
\monend   
the so-called ``dual entropy flux'' $ \zeta^*  (\varphi) $. It is  
defined with the help of the ``physical flux'' $\, F(\smb) \,$ 
according to  
\moneqstar
   \zeta^* (\varphi) \, = \,   
\varphi \, \smb \, F (W) \,\,-\,\,  \zeta (W)\,,  
 \monendstar  
with the condition   (\ref{deriv-entrop})    
as previously.   Then
$ \,   {\rm d} \zeta^* (\varphi)   =  {\rm d} \varphi \,\smb \,  F (W)  \,$ 
%
and the physical flux $ \, F(W) \, $   is the Jacobian of the  dual entropy flux.
In other terms, all the physics associated with the conservation laws  (\ref{hcl}) 
can be expressed in terms of the dual entropy $ \,  \eta^* \,  $ 
and of the dual entropy flux  $ \,  \zeta^* . $ 
  The example of  Burgers equation  (\ref{abs-01}) 
with the quadratic entropy 
and associated flux  
gives without difficulty 
\moneq 
  \displaystyle \eta^*(\varphi)   =  {{\varphi^2}\over{2}} \,\,, \quad 
   \zeta^*(\varphi)    =    {{\varphi^3}\over{6}} \, .   
\label{dual-burgers}    \monend   
%

\monitem 
Independently of the framework relative to hyperbolic conservation laws, 
the Boltzmann equation with discrete velocities has been studied  
by Broadwell   \cite{Br64} 
  (see also Gatignol  \cite{Ga70} and Cabannes  \cite{Ca80}). 
 In this contribution, we write this model  
for $(J \! + \! 1)$ velocities   in one space dimension~: 
\moneq
 \partial_t f_j + v_j  \, \partial_x  f_j = Q_j(f)  
\,  \, ,   \quad 0 \leq j \leq   J \, . 
\label{lbe}   \monend  
The unknown quantity $ \, f_j (x, \, t) \,$ is
the density of particles at point $x$ and time $t$ 
with a discrete velocity $ \, v_j . \,$  
We have for example  $J \! = \! 2$ for the D1Q3 lattice Boltzmann scheme 
(presented in section 4). 
The equation (\ref{lbe}) admits $N$ 
microscopic collision invariants $ \, M_{kj} \,   $:
\moneqstar 
\sum_j M_{kj} \, Q_j(f)  \,=\, 0   \, , \qquad   1 \leq k \leq N
 \monendstar   
and $N=1$ for a scalar ({\it e.g.} Burgers) equation.  
The $N$ first conserved moments : 
\moneq 
W_k \equiv \sum_j M_{kj} \, f_j  \, ,  \quad   1 \leq k \leq N 
\label{conserv-k}   \monend  
satisfy a system   of conservation laws~:
\moneq
 \partial_t W_k +      \partial_x \,  \Big(
  \sum_j  M_{kj} \, v_j  \,    f_j \Big) = 0    \, ,  \quad   1 \leq k \leq N  . 
\label{edp}    \monend   
Of course, we make the hypothesis that 
 this system admits  a mathematical entropy $ \, \eta(W) \,$
with an associated entropy flux $ \, \zeta (W) . \,$ 
We denote by $ \, \varphi \, $ the derivative of the entropy
({\it id est} $ \, {\rm d} \eta =  \varphi \,\smb \, {\rm d}W$)
%
and by $\, M_j \in \R^N  \,$ the vector of components $ \, M_{kj} \,$
(with $k$ running from $1$ to $N$). 
%
Then the following   scalar expression :
\moneq \label{entroparticle}  
 \varphi \, \smb \, M_j \, \equiv \, 
\sum_{k=1}^N   \, \varphi_k  \,  M_{kj}  \,, \quad 0 \leq j \leq J \,, \, 
 \monend 
is well defined. 
%
In some sense, the vector $ \, \varphi \in \R^N  \,$ can be split into 
$ \, J \! + \! 1  \, $ (with $J \geq  N$) {\bf scalar} 
contributions $ \,  \varphi \, \smb \, M_j \,$ associated with the 
particle distribution  of the Boltzmann method. 
In the following, we denote this contribution as 
%
 the ``$j^{\rm o}$ 
particle component  of the entropy variables''. 


\monitem The link between the   Boltzmann models  and the entropy variables
has been first proposed by  Perthame \cite{Pe90}. We follow here 
the approach developed by Bouchut \cite{Bo03}. 
We say that the ``dual entropy approach'' is satisfied if we 
suppose that there  exists 
$J$   {\bf convex}   scalar functions $ \, h_j^* \,$ such that 
\moneq  
 \sum_j  \,  h_j^* \big(  \varphi \, \smb \, M_j \big)   \,\equiv \, 
\eta^* (\varphi)   \,, \quad  
\sum_j  \, v_j   \,  h_j^* \big(  \varphi \, \smb \, M_j \big)  
 \,\equiv \,  \zeta^* (\varphi)  \,, \quad  \forall \varphi   . 
\label{rep-entrop}    \monend   
We introduce $ \, \, 
 h_j(f_j) \,\equiv \, \sup_{y} \, \big( y \, f_j -  h_j^* (y) \big) \, $ 
  the Legendre dual of the convex function $ \,  h_j^* (\smb) .\,$ 
The function $\, h_j(\smb) \,$ is a real scalar convex function 
and   we can write here the relation (\ref{deriv-entrop}) 
making for each $ \, j \,$ the link between $ \, f_j \,$ and 
$ \,  \varphi \, \smb \, M_j  \,$ 
under the scalar form 
\moneq  
h'_j \big( f_j \big) = \varphi \, \smb \, M_j \,, \quad 0 \leq j \leq J  \, . 
\label{derive-hj}    \monend   
The so-called microscopic entropy 
\moneqstar 
 H(f) \equiv \sum_j  h_j(f_j) 
\monendstar 
 is a convex function in the domain where  the $ h_j$'s are convex. 
When  the hypothesis (\ref{rep-entrop}) is satisfied, 
we can prove  a discrete version of the Boltzmann   H-theorem. If  
\moneq
 \sum_j \, h'_j (f_j) \, Q_j(f) \leq 0 , \,\,
\label{hyp-thH} \monend 
we  have dissipation of the microscopic entropy~: 
\moneq
\partial_t H(f)  + 
\partial_x \,  \Big(  \sum_j   v_j    \, h_j (f_j) \Big) \, \leq \, 0  
\label{th-H}    \monend  
and this function is a natural Lyapunov function. 
The equilibrium distribution   $ \, f_j^{\rm eq} (W) \, $  
is then defined by 
\moneq
  f_j^{\rm eq} (W)  \,\equiv \, \big( h_j^* \big) ' 
 \big(  \varphi \,  \smb \, M_j \big)   \,, \quad 0 \leq j \leq J 
\label{flux-eq}   \monend  
because 
the relation 
(\ref{dual-etat})  holds. 
Then we recover the   Karlin  {\it et al} \cite{KGSB98}  minimization property~:  
$ \,\,   H(f) \geq H(f^{\rm eq}) \,\,$ 
 for each $ \,  f \,$    such  that 
$\,\,  \sum_j \, M_{kj} \,  f_j \,=\, \sum_j \, M_{kj} \,  f_j^{\rm eq} 
\equiv W_k \,$ with $ \, 1 \leq k \leq N . $ 

%
\monitem  
By differentiation of the relations (\ref{rep-entrop}) relative to 
the entropy variable $\, \varphi \,$ and taking into account the previous
relations (\ref{flux-eq}), we have the necessary equilibrium conditions
$ \, \sum_j \, M_j \,   f_j^{\rm eq} =  W \, $ 
and $ \,  \sum_j \, v_j \,  M_j \,   f_j^{\rm eq}  = F(W) $. 
In other terms,   
%
the conserved variables  
are given  by the relations (\ref{flux-eq})(\ref{conserv-k}) 
and the macroscopic fluxes by 
%
\moneqstar
  F_k(W) \,\equiv 
\,\sum_j \, M_{kj} \, v_j \,   f_j^{\rm eq}
 \, \,, \,\,\, 1 \leq k \leq N \, . 
\monendstar
The macroscopic entropy and associated entropy fluxes   
satisfy 
\moneqstar
  \eta(W) \, = \, \sum_j  h_j \big( 
  f_j^{\rm eq} \big) \,, \quad    \zeta (W) \, = \, \sum_j  v_j  \, 
h_j \big(    f_j^{\rm eq} \big)  \, . 
\monendstar
When the  Boltzmann equation with discrete velocities satisfies  the so-called 
BGK hypothesis \cite{BGK}, {\it id est} 
\moneq 
 Q_j(f) =   {{1}\over{\tau}} \big(  f_j^{\rm eq} - f_j \big)   \,, \quad 
0 \leq j \leq J 
\label{bgk-hyp} \monend
for some constant  $  \tau > 0 , \, $  
the Boltzmann H-theorem is satisfied. We give the proof for 
completeness~: 
we first have  the following convexity inequality 
\moneqstar  
 \Big(   h_j' \big(f_j^{\rm eq} \big)  -   h_j' \big(f_j \big) 
\Big) \, \, \Big(  f_j ^{\rm eq} - f_j  \Big) \, \geq \, 0 \,,    \quad 
0 \leq j \leq J \, . 
\monendstar
If the BGK hypothesis (\ref{bgk-hyp}) occurs, we have 
 by summation over $j$, 
 
\smallskip  \noindent
$ \displaystyle \tau \, \sum_j h_j' \big(f_j \big) \,    Q_j(f)  
\, = \,  \sum_j  h_j' \big(f_j \big) \,  \Big(  f_j ^{\rm eq} - f_j  \Big) 
\, \leq \,  \sum_j  h_j' \big(f_j^{\rm eq} \big) \,  
 \, \Big(  f_j ^{\rm eq} - f_j     \Big) 
 \, =  \,   $   

\hfill 
$ \displaystyle   \, =  \,  \sum_j   \big(   \varphi \,  \smb \, M_j   \big) \, 
\Big(  f_j ^{\rm eq} - f_j     \Big)  \, =  \, 
 \varphi \,  \smb \,  \sum_j     M_j    \, 
\Big(  f_j ^{\rm eq} - f_j     \Big) \, = \, 0  \,   $   

\smallskip \noindent 
and due to (\ref{derive-hj}), the hypothesis (\ref{hyp-thH}) 
is satisfied. In consequence  the H-theorem is established in this case. 

 
\monitem As a summary of this mathematical section, 
we  explicit the dual entropy approach
 in the case of the Burgers
equation (\ref{abs-01}) equipped with a quadratic entropy. 
If 
there  exists {\bf convex} functions 
$ \, h_j^*(\varphi) \,$ of the entropy variable $ \, \varphi \,$  such that 
\moneq   \label{abs-03}   
   \sum_j  h_j^*(\varphi)   \, \equiv  \,  \eta^* (\varphi)  \, = \,  
  {{\varphi^2}\over{2}} \,, \quad   
 \,   \sum_j  v_j \, h_j^*(\varphi)    \, \equiv  \, \zeta^* (\varphi) \, = \,  
  {{\varphi^3}\over{6}}   \monend
then the equilibrium 
$\, f_j^{\rm eq}(u) \equiv  {{{\rm d} h_j^*}\over{{\rm d} \varphi}}  \,$  
defines a {\bf stable} approximation in a sense detailed in 
Chen {\it et al} \cite{CLL94} and extended by   Bouchut \cite{Bo99, Bo04}.

\bigskip \bigskip  \noindent {\bf \large 3) \quad Particle  decompositions for the Burgers
  equation}

\quad  
We propose in this contribution to construct  kinetic decompositions
of a scalar variable in order to solve the Burgers equation in cases where weak solutions
can occur, {\it id est} when shock waves can be developed. We consider only the simple
D1Q3 stencil with three discrete velocities $-\lambda$, $0$ and $\lambda$.
Recall that the scalar  $\lambda \equiv {{\Delta x}\over{\Delta t}}$ is 
a fundamental numerical parameter that is very often taken equal to  unity 
by lattice Boltzmann scheme users (see {\it e.g.} \cite {LL00}). 
For   the Burgers equation  (\ref{abs-01})
a  possible mathematical entropy is the quadratic one   (\ref{quad-entrop}).
The dual entropy $   \,   \eta^*(\varphi) \,$ and the associated dual entropy flux 
$ \,   \zeta^*(\varphi) \,      $ 
are given according to the relations  (\ref{dual-burgers}).  
Due to the 
framework of dual entropy approach
 proposed in the previous section, we search 
three  convex functions $ \, h_+^* (\varphi)   \,, $  
    $ \, h_0^* (\varphi)   \, $ and $ \, h_-^* (\varphi)   \,$
such that  (\ref{abs-03}) holds, {\it id est} for D1Q3~:  
\moneq   \label{conds-burgers}    
   h_+^* (\varphi) +   h_0^* (\varphi)  +   h_-^* (\varphi) 
\equiv    {{\varphi^2}\over{2}} \,, \quad  
 \lambda \, \big(  h_+^* (\varphi) \,-\,   h_-^* (\varphi) \big) 
\equiv     {{\varphi^3}\over{6}}.    \monend
%

\monitem 
A first possible solution of the previous system 
consists in introducing some  parameter  $\alpha$ such that
 $   0 < \alpha \leq 1  . \, $  Then we  consider the particular function 
\moneq   \label{h-zero-alpha}  
 h_0^* (\varphi) \,=\, (1-\alpha) \, {{\varphi^2}\over{2}} \, . 
  \monend
%
Of course, if $ \, \alpha = 1 ,\,$ this  function $ \, h_0^* (\smb) \,$
is singular. In this case, we switch from D1Q3 to D1Q2 scheme, as presented  
in the following of this contribution. 
%
Due to (\ref{conds-burgers}), the two other dual functions 
$\,  h_+^* (\varphi)   \, $ and $ \, h_-^* (\varphi)   \,$ 
are determined~: 
\moneq   \label{h-pm-alpha}  
 h_+^* \,=\, \alpha \, {{\varphi^2}\over{4}} 
\,+\,  {{\varphi^3}\over{12 \, \lambda}} \,,  \quad 
 h_-^* \,=\, \alpha \, {{\varphi^2}\over{4}} 
\,-\,  {{\varphi^3}\over{12 \, \lambda}}  \, . \monend
The associated dual functions can be written explicitly 
 without particular difficulty~:
\moneq   \label{h-de-f} 
\left\{ \begin{array} [c]{cl}  
\displaystyle   h_+ (f_+) \, & = \,\,  \displaystyle  
{{\lambda^2}\over{6}} \, \Big[ 
\Big( \alpha^2 + 4 \,{{f_+}\over{\lambda}} \Big)^{3/2} \,- \, 
6 \, \alpha \, {{f_+}\over{\lambda}}   \,- \,   \alpha^3    \Big]    
\\ \displaystyle   h_0 \big( f_0 \big)  \, & = \,\, \displaystyle 
   {{1}\over{2 \, (1 - \alpha)}}  \, f_0 ^2 
\\ \displaystyle    h_- \big( f_- \big)  \, & = \,\,   \displaystyle 
{{\lambda^2}\over{6}} \, \Big[  
\Big( \alpha^2 - 4 \,{{f_-}\over{\lambda}} \Big)^{3/2} 
+ \, 6 \, \alpha \, {{f_-}\over{\lambda}} - \alpha^3   \Big] \, .  
\end{array} \right.     \monend
The three functions $\, h^*_j \,$ introduced in  (\ref{h-zero-alpha})
and  (\ref{h-pm-alpha}) are   convex   when  
\moneq   \label{convex-alpha}  
 \  \mid \! \varphi \!\mid \, \leq \, \alpha \, \lambda   
\monend
and the relation (\ref{convex-alpha}) can be interpreted as a 
Courant-Friedrichs-Lewy stability condition~: 
\moneqstar
 \Delta t \, \leq \, 
{{\alpha}\over{\mid \! u \!\mid}}   \, \Delta x \, . 
\monendstar
%
The dual entropy approach contains in particular the  numerical
stability condition (\ref{convex-alpha}). 
%
The stability is in fact defined as the domain of convexity 
 of the dual functions $ \, h_j^* $ presented algebraically by  relations 
(\ref{h-zero-alpha}) (\ref{h-pm-alpha}) and illustrated 
in  Figure~1. 
The explicit  determination of the  equilibrium distribution is then a 
consequence of the relation (\ref{flux-eq}) taking also into account that 
 $ \, \varphi \equiv u \, $ for the quadratic entropy. We have  
\moneq \label{centre-d1q3-f}
   f_+^{\rm eq} (u)  = 
{{\alpha}\over{2}} \, u  +  {{u^2}\over{4\, \lambda}}  \,, \quad 
  f_0^{\rm eq}  =   (1 - \alpha) \, u   \,, \quad  
 f_-^{\rm eq}   =  {{\alpha}\over{2}} \, u  - {{u^2}\over{4\, \lambda}}  \, . 
\monend
%
%

\smallskip 
\centerline { \includegraphics[width=.50 \textwidth, angle=-90]
{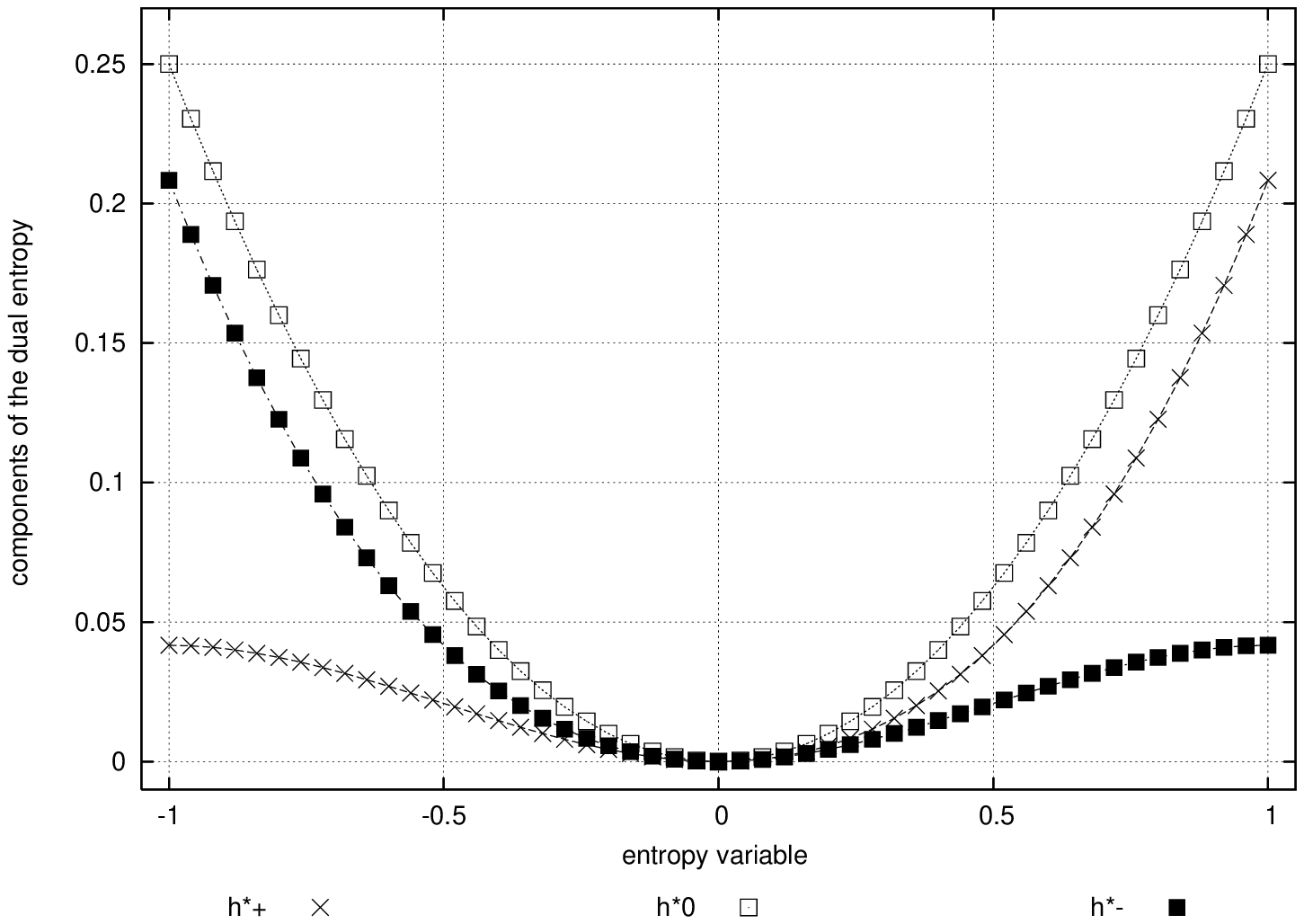}}  

\smallskip \noindent  {\bf Figure 1}. \quad 
Kinetic decomposition (\ref{h-zero-alpha}) (\ref{h-pm-alpha}) 
of the dual entropy   for the Burgers 
equation  with   a ``centered''     D1Q3 scheme ($\alpha = {{1}\over{2}}$).   
\smallskip 

\smallskip 
\centerline { \includegraphics[width=.50    \textwidth, angle=-90]
{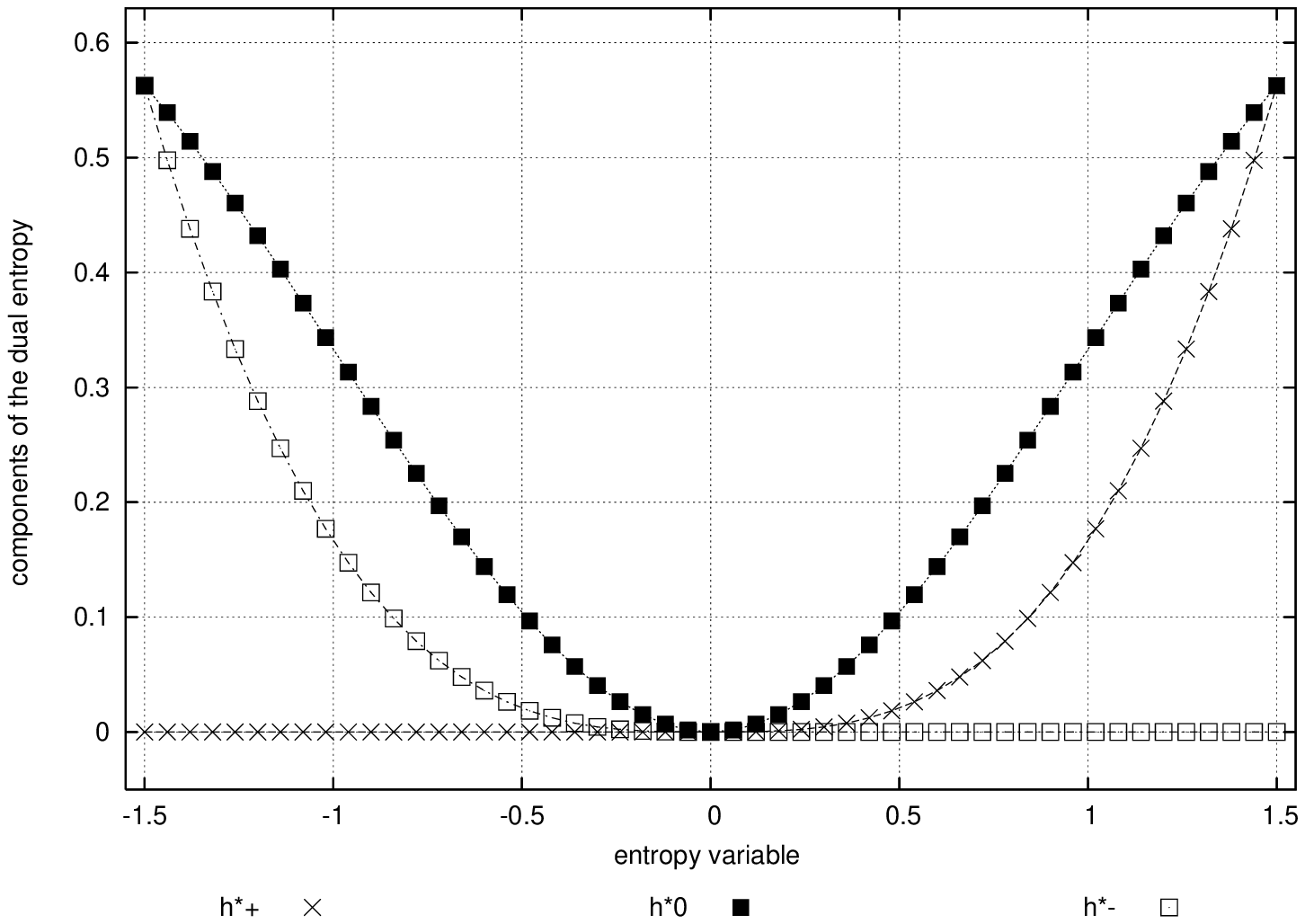}}   

\smallskip \noindent  {\bf Figure 2}. \quad  
Kinetic decomposition for Burgers equation, 
equilibria (\ref{upwind}) for the lattice Boltzmann upwind scheme D1Q3.   
\smallskip 


\monitem 
Another solution of the previous system  (\ref{conds-burgers}) can be obtained as
follows.  
Derive the two  relations in  (\ref{conds-burgers}) two times. Then  
\moneqstar  
  \big(h_+^*\big)''(\varphi)  \,=\,  \big(h_-^*\big)''(\varphi)  \,+ \, 
  {{\varphi}\over{\lambda}}  \, , \quad 
 \big(h_0^*\big)''(\varphi)    \,+\, 2 \,  \big(h_-^*\big)''(\varphi)    \,=\, 
1 -    {{\varphi}\over{\lambda}} \, . 
\monendstar
In order to have a better stability property than the condition 
 (\ref{convex-alpha}) obtained previously, 
we try to enforce the convexity condition  
$ \,   \big(h_j^*\big)''(\varphi) \geq 0  \,$ if 
$\,   \mid \! \varphi \!\mid \, \leq \,  \lambda  \, $
instead of (\ref{convex-alpha}).  
For     $\, \varphi \leq 0 ,\,$ we propose to replace the inequality  
$ \,  \big( h_+^*\big)''(\varphi) \equiv   \big(h_-^*\big)''(\varphi)  
+   {{\varphi}\over{\lambda}} \,\geq \, 0 \,$
 by an equality. Then  
$ \,   \big(h_-^*\big)''(\varphi) =  - {{\varphi}\over{\lambda}}   \,$ 
if $ \,    \varphi \leq 0 . \, $ We   deduce 
$ \,  \big(h_+^*\big)'' (\varphi) =  0 \,  $ and  
$ \, \big(h_0^*\big)''  (\varphi) =  1 +  {{\varphi}\over{\lambda}} \, $ 
if $ \,   \varphi \leq 0 . \, $  
With analogous arguments, we obtain 
$ \,  \big(h_+^*\big)'' (\varphi) =  {{\varphi}\over{\lambda}}  , \, $ 
$ \,  \big(h_0^*\big)'' (\varphi) =   1 -   {{\varphi}\over{\lambda}}    \, $
and 
$ \,   \big(h_-^*\big)'' (\varphi) = 0  \, $ 
when  $ \, \varphi \geq \lambda . \,$
%
We   construct in this way an ``upwind''  distribution 
for the decomposition of the dual entropy:   
\begin{equation} \label{upwind}   
h_+^* (\varphi)  = 
\left\{  \begin{array} [c] {l} \displaystyle   {{\varphi^3}\over{6 \, \lambda}} 
      \\ ~ \\ \displaystyle  0  \end{array} \right. , \quad  h_0^* (\varphi) = 
\left\{ \begin{array} [c]  {l}  \displaystyle 
  {{\varphi^2}\over{2}}  -  {{\varphi^3}\over{6 \, \lambda}} 
  \\ ~  \\  \displaystyle   {{\varphi^2}\over{2}}  +  {{\varphi^3}\over{6 \, \lambda}} 
 \end{array} \right.   , \quad h_-^* (\varphi) = 
\left\{ \begin{array} [c]{cl}  \displaystyle   0 , 
   &  \quad  \varphi \geq 0   \\ ~&~  \\   \displaystyle 
-  {{\varphi^3}\over{6 \, \lambda}} , &   \quad  \varphi \leq 0 .  \end{array} \right. 
\end{equation}  
It is presented in  Figure~2. 
The associated  equilibrium distribution (\ref{flux-eq}) 
 takes the form 
\begin{equation} \label{upwind-f}   f_+^{\rm eq} (u)  = 
\left\{  \begin{array} [c] {l} \displaystyle  {{u^2}\over{2 \, \lambda}}  
      \\ ~ \\ \displaystyle  0  \end{array} \right. , 
\quad  f_0^{\rm eq}  (u)  = 
\left\{ \begin{array} [c]  {l}  \displaystyle 
 u -  {{u^2}\over{2 \, \lambda}} 
  \\ ~  \\  \displaystyle u +  {{u^2}\over{2 \, \lambda}} 
 \end{array} \right.   
, \quad   f_-^{\rm eq}  (u)  = 
\left\{ \begin{array} [c]{cl}  \displaystyle   0 , 
   &  \quad  u \geq 0   \\ ~&~  \\   \displaystyle 
-  {{u^2}\over{2 \, \lambda}}  , &   \quad  u  \leq 0  \, .  \end{array} \right. 
\end{equation} 
By considering the Legendre duals of the relations (\ref{upwind}), we have 
\moneq   \label{h-engquist-de-f} 
\left\{ \begin{array} [c]{cll}  
\displaystyle   h_+ (f_+) \, & = \,\,\,\,   \displaystyle  
{{2}\over{3}} \, f_+ \, \sqrt{2 \, \lambda \, f_+} &  {\rm with} \, f_+ \geq 0   
\\ \displaystyle   h_0 \big( f_0 \big)  \, & = \, \quad  
 \displaystyle {{\lambda^2}\over{3}} \,   
\Big[ \Big( 1 - 2 \,{{ \mid \!\! f_0\!\! \mid }\over{\lambda}}  \Big)^{3/2} 
 + 3\, {{ \mid \!\! f_0\!\! \mid }\over{\lambda}} - 1  \Big] 
\quad &    {\rm with} \, f_0 \in \R    
\\ \displaystyle    h_- \big( f_- \big)  \, & = \,    \displaystyle 
-{{2}\over{3}} \, f_- \, 
\sqrt{-2 \, \lambda \, f_-} &  {\rm with} \, f_- \leq 0   \, .  
\end{array} \right.     \monend
%


\monitem  
We observe that if $ \, \alpha = 1 \,$ for the ``centered'' equilibrium for D1Q3
Burgers scheme, the null velocity does not contribute to the equilibrium 
because $ \, h_0(\varphi) \equiv 0 \,; $ this vertex of null velocity is no more
active. In that case, we obtain a  D1Q2 centered  lattice Boltzmann scheme for Burgers
equation. Then 
\moneq  \label{hstar-d1q2}    
h_+^* (\varphi)  =  {{\varphi^2}\over{4}} 
\,+\,  {{\varphi^3}\over{12 \, \lambda}} \,,  \quad 
 h_-^* \,=\,  {{\varphi^2}\over{4}} 
\,-\,  {{\varphi^3}\over{12 \, \lambda}}   \, . 
\monend
These two   functions represented in Figure~3 
are   convex    if 
\moneq  \label{stabi-d1q2}    
 \mid \! \varphi \!\mid \, \leq \,  \lambda    
\monend
and the associated  
Courant-Friedrichs-Lewy     stability condition   states as follows 
\moneqstar
  \Delta t \, \leq \, 
{{1}\over{\mid \! u \!\mid}}   \, \Delta x  \, . 
\monendstar 
The dual equilibrium entropy function defined at relations 
(\ref{hstar-d1q2}) are represented on Figure~3. 
The associated components $ \, h_+(f_+) \,$ and  $ \, h_-(f_-) \,$
of the microscopic entropy follow from  (\ref{h-de-f}) 
in  the particular case $\, \alpha = 1 .$ Observe that 
 $ \, h_0(f_0) \,$ is no more defined which is coherent 
with a choice of a ``D1Q2'' lattice Boltzmann scheme. 
The associated equilibrium particle  distribution 
 is obtained according to 
\moneq \label{distri-equil-d1q2}
   f_+^{\rm eq} (u) =  {{1}\over{2}} \, u  + {{u^2}\over{4\, \lambda}} 
\,, \quad 
  f_-^{\rm eq} = {{1}\over{2}} \, u  - {{u^2}\over{4\, \lambda}}  \, . 
  \monend  
%
\smallskip  
\centerline { \includegraphics[width=.50    \textwidth, angle=-90 ]
{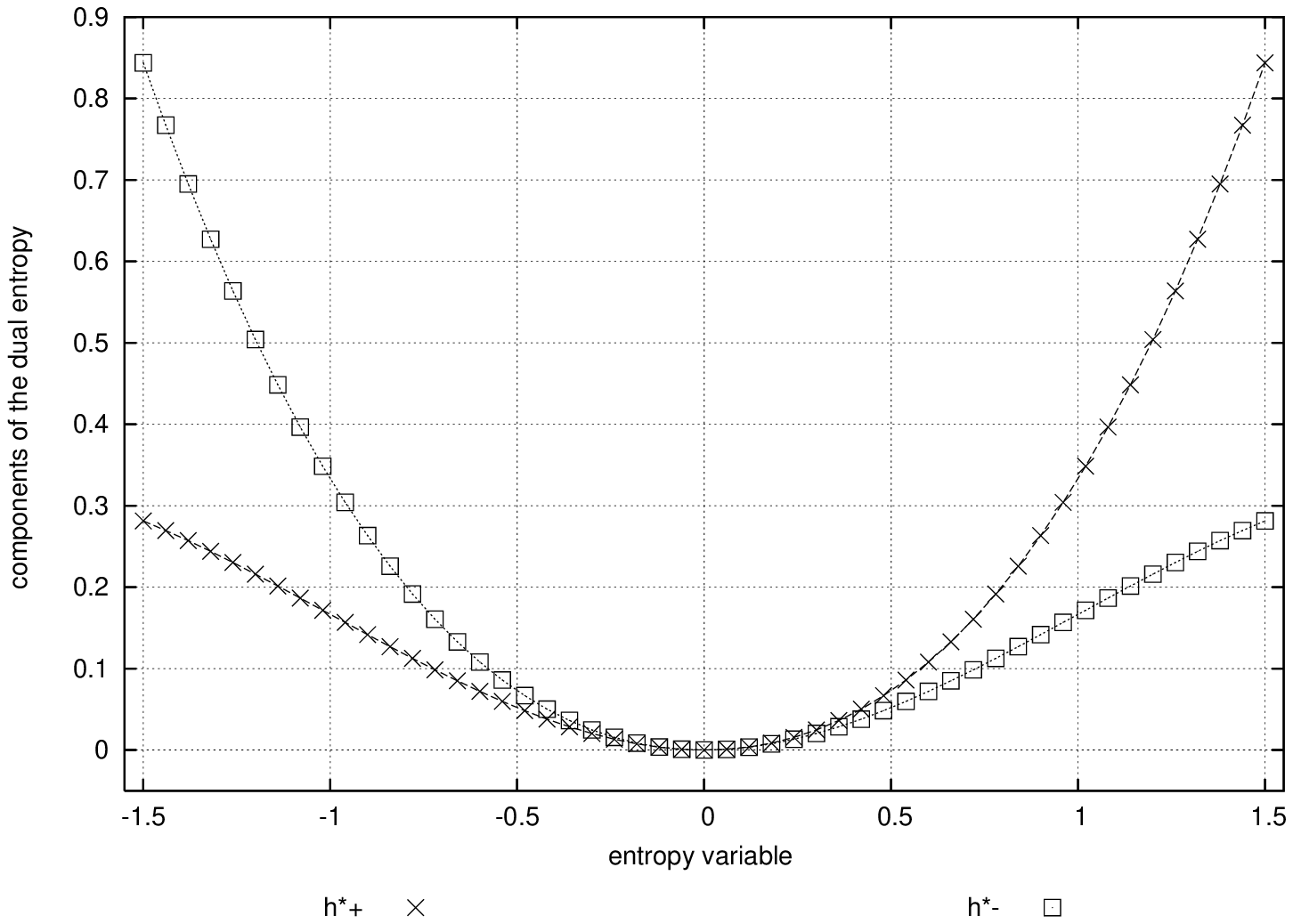}}    

\smallskip \noindent  {\bf Figure 3}. \quad  
Kinetic decomposition for Burgers~;  D1Q2   centered scheme.   
\smallskip 
%

\bigskip \bigskip  \noindent {\bf \large 4) \quad D1Q3 lattice Boltzmann scheme }      
 
\quad  
As developed in the preceding section, we here consider  three  examples of stable  equilibria 
in the context of the lattice Boltzmann scheme. 
More precisely, following the  approach proposed by   d'Humi\`eres \cite{DDH92}, 
we discretize in space and time the 
Boltzmann equation  with discrete velocities  (\ref{lbe})   in the following way.
We introduce a matrix $M$ that links  particle densities 
$ \, f_j \,$ $(j=- ,\, 0 ,\, +)$  and moments~$\, m_k$. 
For the simple D1Q3 lattice Boltzmann scheme, we obtain 
\moneq  \label{abs-04}   
m \, \equiv \,  M \, \smb \, f  \,, 
\quad  M \,=\, \left( \begin{matrix} 1 & 1 & 1\\ -\lambda & 0 & \lambda  \\ 
 \lambda^2 & 0 &  \lambda^2   \end{matrix}\right)  \,,  
\quad  u \equiv f_{-1} + f_0 + f_1 = m_1 \,.  
\monend

\monitem 
The first equilibrium  (\ref{centre-d1q3-f}) can be translated  in terms of moments 
under the form 
\moneqstar 
m^{{\rm eq},1} \, \equiv  \, \Big( u \,,\,  {{u^2}\over{2}} \,,\, 
\alpha \, \lambda^2  \, u   \Big) ^{\rm \displaystyle t} \, . 
\monendstar
When using the ``upwind'' equilibrium (\ref{upwind-f}), we obtain an other 
possible value for moments at equilibrium~:  
\moneqstar 
m^{{\rm eq},2} \, \equiv \, \Big( u \,,\,  {{u^2}\over{2}} \,,\, \lambda \, 
{\rm sgn}(u) \, {{u^2}\over{2}} \Big) ^{\rm \displaystyle t}  \, . 
\end {equation*} 
The simpler scheme D1Q2 corresponds to the first equilibrium  (\ref{centre-d1q3-f}) 
with the particular value $ \alpha = 1 \,$ as proposed in relations 
 (\ref{distri-equil-d1q2}). We have only two components in this case~:
\moneqstar
m^{{\rm eq},3} \, \equiv \, \Big( u \,,\,  {{u^2}\over{2}}  \Big) ^{\rm \displaystyle t}  \, . 
\monendstar

\monitem  
The relaxation step is nonlinear and local in space~: 
\moneq  \label{relax-moments}
 m^*_1  =  m^{\rm eq}_1 = u \, , \quad  
 m^*_k  = m_k + s_k \, ( m^{\rm eq}_k - m_k) \, \, {\rm for} 
\,  \, k \geq 2 ,\, 
\monend 
 with $ \, s_2 = s_3 = 1.7 \,$ 
in our simulations unless otherwise stated. 
%
For  nonlinear hyperbolic systems (\ref{hcl}) 
of two conservation laws in one space dimension, 
the moments $\, m_1 \,$ and  $\, m_2 \,$
are at equilibrium and the relation (\ref{relax-moments}) is written in this case 
\moneq  \label{relax-moments-system}
 m^*_1  =  m^{\rm eq}_1 = W_1 \, , \quad  
 m^*_2  =  m^{\rm eq}_2 = W_2 \, , \quad  
 m^*_3  = m_3 + s_3 \, ( m^{\rm eq}_3 - m_3) \, . 
\monend 
%
The particle distribution  $ \, f^*_j \,$  after relaxation 
is obtained by inversion of  relation   (\ref{abs-04})~: $ \, f^* = M^{-1}  \smb \ m^*.$
%
 The time iteration of the 
scheme follows the characteristic directions of velocity $ \, v_j $ : 
\moneqstar 
f_j(x ,\, t+ \Delta t) =  f_j^*(x-v_j \, \Delta t ,\, t)  \, .  
\monendstar
This advection step is linear and associates the node $x$ with its  neighbors.

\monitem 
In \cite{Du08}  we have observed that a one-dimensional lattice Boltzmann scheme can be
interpreted with the help of finite volumes. 
In the case considered here, we have 
\moneqstar 
{{1}\over{\Delta t}} \Big( u(x,\, t+\Delta t) -   u(x,\, t) \Big) 
+ {{1}\over{\Delta x}} \Big[ \psi \Big(  x+ {{\Delta x}\over{2}}, \, t \Big) 
-  \psi \Big(  x - {{\Delta x}\over{2}}, \, t \Big) \Big] = 0  
\monendstar
with a numerical flux $ \,  \psi \big(  x+ {{\Delta x}\over{2}}, \, t \big) \,$
at the interface between the vertices $x$ and $ \, x + \Delta x \,$ 
defined according to  
\moneq \label{flux-psi} 
\psi \Big(  x+ {{\Delta x}\over{2}}, \, t \Big) = \lambda \, \Big( 
f_+^* (x, \, t) - f_-^* (x + \Delta x, \, t)  \Big) \, .  
\monend 
We observe that the resulting lattice Boltzmann scheme is {\bf not}  a traditional finite volume
scheme (in the sense proposed {\it e.g.} in \cite{DD05})  
if $ \, (s_2,\, s_3) \not= (1,\, 1) \,$ because the distribution 
of particles after collision $ \, f^* \,$ is also a function of the two (or one in the
D1Q2 scheme) other nonconserved moments $ \, m_2 \,$ and $ \, m_3  \,$  
as described in  relations (\ref{relax-moments}). 
%
On the contrary, the lattice Boltzmann method is mainly 
 a particle method with given velocities, as 
analyzed  {\it e.g.} in   Junk {\it al.} \cite {JKL05} with an
asymptotic expansion technique.  
%
Nevertheless, if $ \,  s_2 =  s_3 = 1 ,\, $ we can give an interpretation of the
associated flux (\ref{flux-psi}) because in this case, $ \, f^*_j \equiv  f^{\rm eq}_j \,$
for all $j$.  

\monitem 
We observe that we can also decompose the ``physical'' flux $\, F(\smb)\,$ 
(see the relation  (\ref{abs-01}) or (\ref{hcl}) in all generality) 
under the form $ \,  F(u) \equiv   F_+(u) + F_-(u)   \, $ 
%
%
with 
\moneq  \label{decomp-flux} 
 F_+(u) = \lambda \, f_+^{\rm eq}(u) \,, \quad 
 F_-(u) = - \lambda \, f_-^{\rm eq}(u) \, . 
\monend 
We have  
$ \,\,  F_+ \big( u(x,\, t) \big) +  F_- \big( u(x + \Delta x,\, t) \big) = 
\lambda \, \big( f_+^{\rm eq}  \big( u(x,\, t) \big) - 
 f_+^{\rm eq}  \big( u(x+ \Delta x ,\, t) \big) \big) \, \,$  
and when  $ \,  s_2 =  s_3 = 1 \, $ the numerical flux $\psi$ introduced in 
(\ref{flux-psi}) admits the classical so-called flux splitting form~: 
\moneq  \label{lien-flux-particules}  
\psi \Big(  x+ {{\Delta x}\over{2}}, \, t \Big) = 
 F_+ \big( u(x,\, t) \big) +  F_- \big( u(x + \Delta x,\, t) \big)  \, . 
\monend
With  this above link between fluxes and particle distributions  (\ref{lien-flux-particules})   
it is  natural to re-interpret,  
with  classical flux decompositions  as  (\ref{decomp-flux}), 
those  proposed in this contribution at relations  
  (\ref{centre-d1q3-f}),   (\ref{upwind-f})  and   (\ref{distri-equil-d1q2}). 
As remarked by Bouchut \cite {Bo10}, the relations
  (\ref{centre-d1q3-f}) and   (\ref{distri-equil-d1q2}) are associated with two  
variants of the Lax-Friedrichs scheme (see {\it e.g.} Lax \cite {La06}) 
whereas the upwind scheme   (\ref{upwind-f}) corresponds  exactly to the 
 Engquist-Osher \cite{EO80} scheme~!

\centerline { \includegraphics[width=.65 \textwidth, angle=-90]
{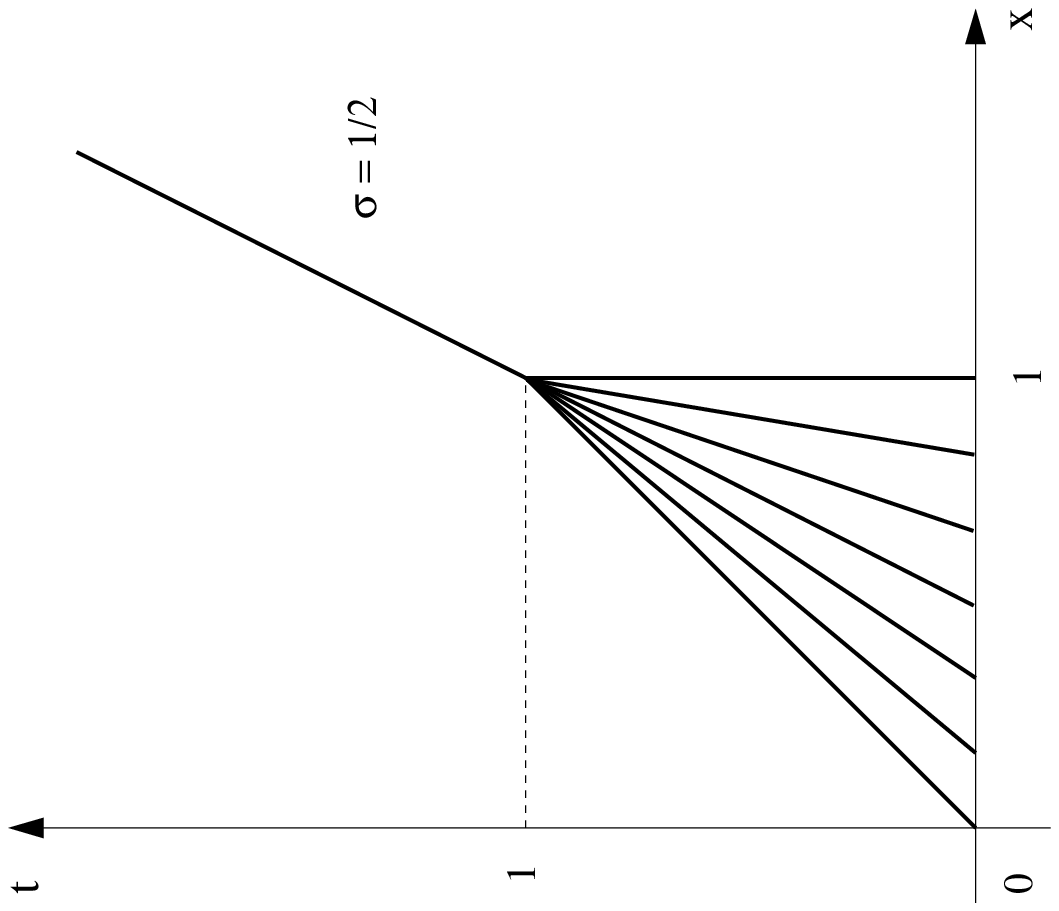} } 

\smallskip \noindent  {\bf Figure 4}. \quad 
  A converging shock wave for the Burgers equation. The decreasing profile 
 (\ref{profil-decroi}) at $t=0$ leads  to an admissible discontinuity at $t=1$. 
Then a shock wave with  velocity $ \, \sigma = {1\over2} \,$ develops. 

\bigskip \bigskip  \noindent {\bf \large 5) \quad Test cases for Burgers nonlinear waves}  

\quad 
We   test the previous numerical schemes for two classical problems~:
a   converging shock wave  and the Riemann problem. We use  the three variants 
(\ref{centre-d1q3-f}), (\ref{upwind-f}) and  (\ref{distri-equil-d1q2}) 
of the lattice Boltzmann scheme for each   problem.

\monitem The first test case concerns a  converging shock wave 
and is displayed in Figure~4. 
At time $t=0$ the initial profile $ \, u_0(x) \,$ is given according to 
\moneq \label{profil-decroi}    u_0(x)  = 
\left\{  \begin{array} [c] {cl}  \displaystyle  1 \quad & {\rm if}  \quad   x \leq 0 \\
1 - x  \quad  & {\rm if} \quad  0 \leq x \leq 1 \\
0  \quad  & {\rm if} \quad   x \geq 1  \, .  \end{array} \right.   
\monend 
%
When $t < 1$ the profile $\, u(x,\, t) \,$ remains a continuous function of space $x$ but
when $t > 1 $ a shock wave with  velocity  $ \, \sigma = {1\over2} \,$ 
is present (see {\it e.g.}  
 \cite{GR96},   \cite{DD05} or \cite{La06}). It is a challenge if a lattice
Boltzmann scheme is able to capture in a systematic way 
such a discontinuous solution. 


\monitem The first experiment (see Figure~5) 
concerns the 
first centered scheme  (\ref{centre-d1q3-f}) and the choice 
 $ \alpha = {1\over2} $   and  $  \lambda = 1.8 $
for the numerical parameters. The result is catastrophic, as depicted on  
Figure~5. 
The scheme is unstable and diverges within a  
very little time after the solution becomes discontinuous. 
The reason is simple {\it a posteriori}.  
Observe that for the previous test case $ \, \alpha = {1\over2} \, $ 
and  particular values 
$ u(x,\,t) \geq 1 $ have to be considered.  But the  
convexity-stability condition  (\ref{convex-alpha}) reads  as   
 $ \, \mid \! u \! \mid  \, \leq \,     {{\lambda}\over{2}} \, $ and 
is incompatible with the chosen numerical values because we take 
$\lambda = 1.8 \,\, $ in the numerical simulation. 
We observe that under   conditions that violate the inequality (\ref{convex-alpha}), 
the lattice Boltzmann scheme is unstable in this strongly nonlinear
situation, even if  we   respect  the linear stability  condition 
\moneq         \label{henon}   
 0 < s_j < 2   
\monend  
proposed initially by   H\'enon \cite{He87}.

\smallskip   
\centerline { \includegraphics[width=.50 \textwidth, angle=-90  ]
{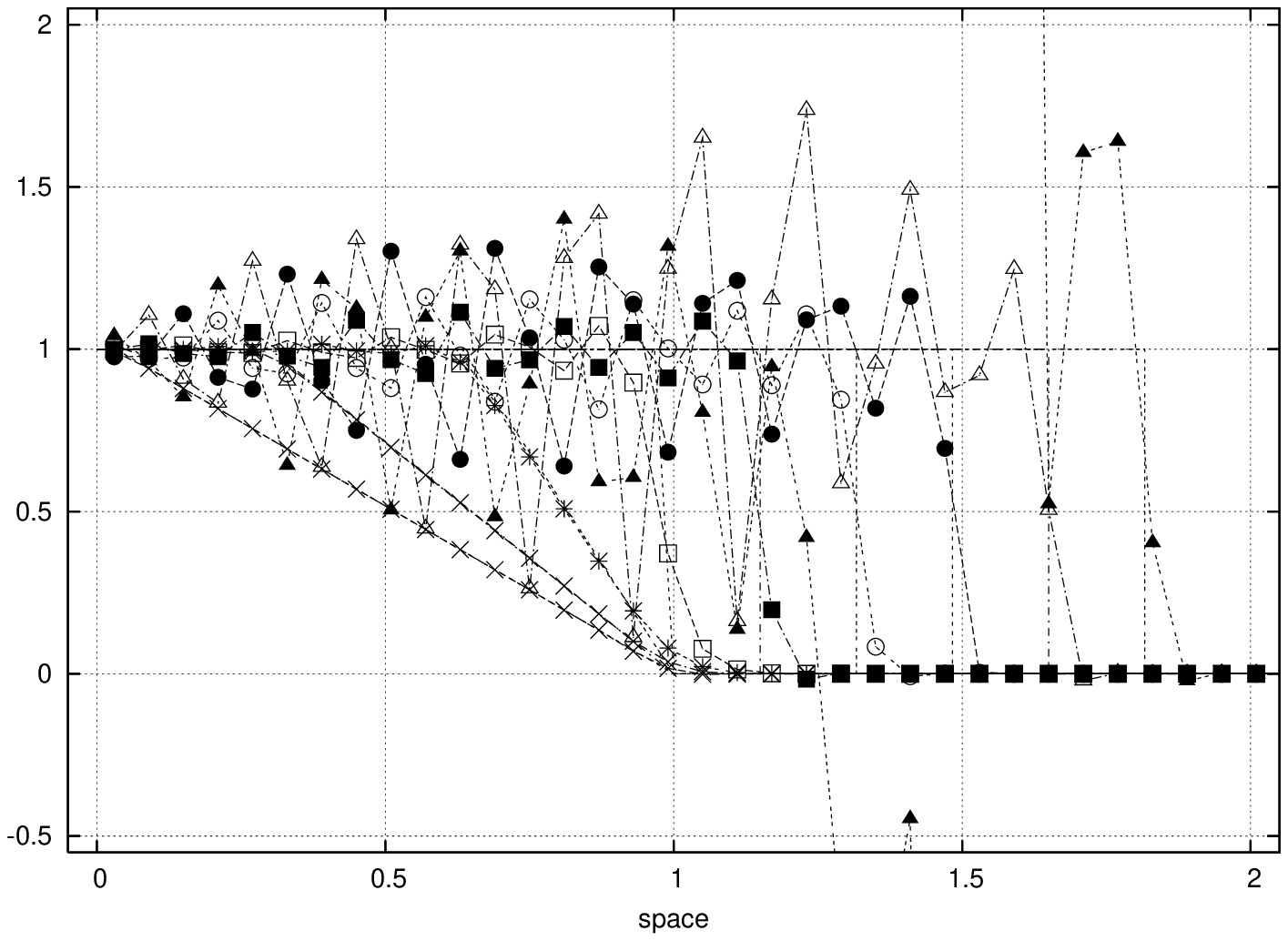} } 
 \smallskip   \smallskip \noindent  {\bf Figure 5}. \quad  
Burgers equation. Instable   D1Q3  lattice Boltzmann simulation for 
a converging shock with equilibrium (\ref{centre-d1q3-f}) 
associated to the parameters $ \alpha = {1\over2} $,  $ s_2 = s_3 = 1.7$ 
and  $  \lambda = 1.8 $. Computed values are displayed  every 10 time steps.  
 \smallskip   \smallskip  
%

%
 
\monitem We repeat the same numerical experiment with 
 a smaller time step. We   take   $ \, \lambda = 3 \,  $
in a second experiment. 
 The condition    (\ref{convex-alpha}) 
is now satisfied and the scheme is stable. The results 
are correct and are presented in Figure~6. 
The shock is spread on 4 to 5 mesh points and
we observe simply an overshoot at the location of the shock wave. 
With the extreme set of  values $ \,  s_2 = s_3 = 2 $
(if we refer to   relation (\ref{henon})),  
the numerical experiment does not give correct results
 because no entropy is dissipated. But the scheme remains stable; 
%
%
the numerical values
remain inside an interval $ \, [-0.4 ,\, 1.7 ] \, $ relatively  close to the set 
$\, [0,\, 1]\,$ of correct values for this particular problem. 
The nonlinear stability condition enters into competition with the 
linear stability condition~(\ref{henon}). 

 \smallskip     \bigskip     \bigskip   
\centerline { \includegraphics[width=.50 \textwidth, angle=-90  ] 
{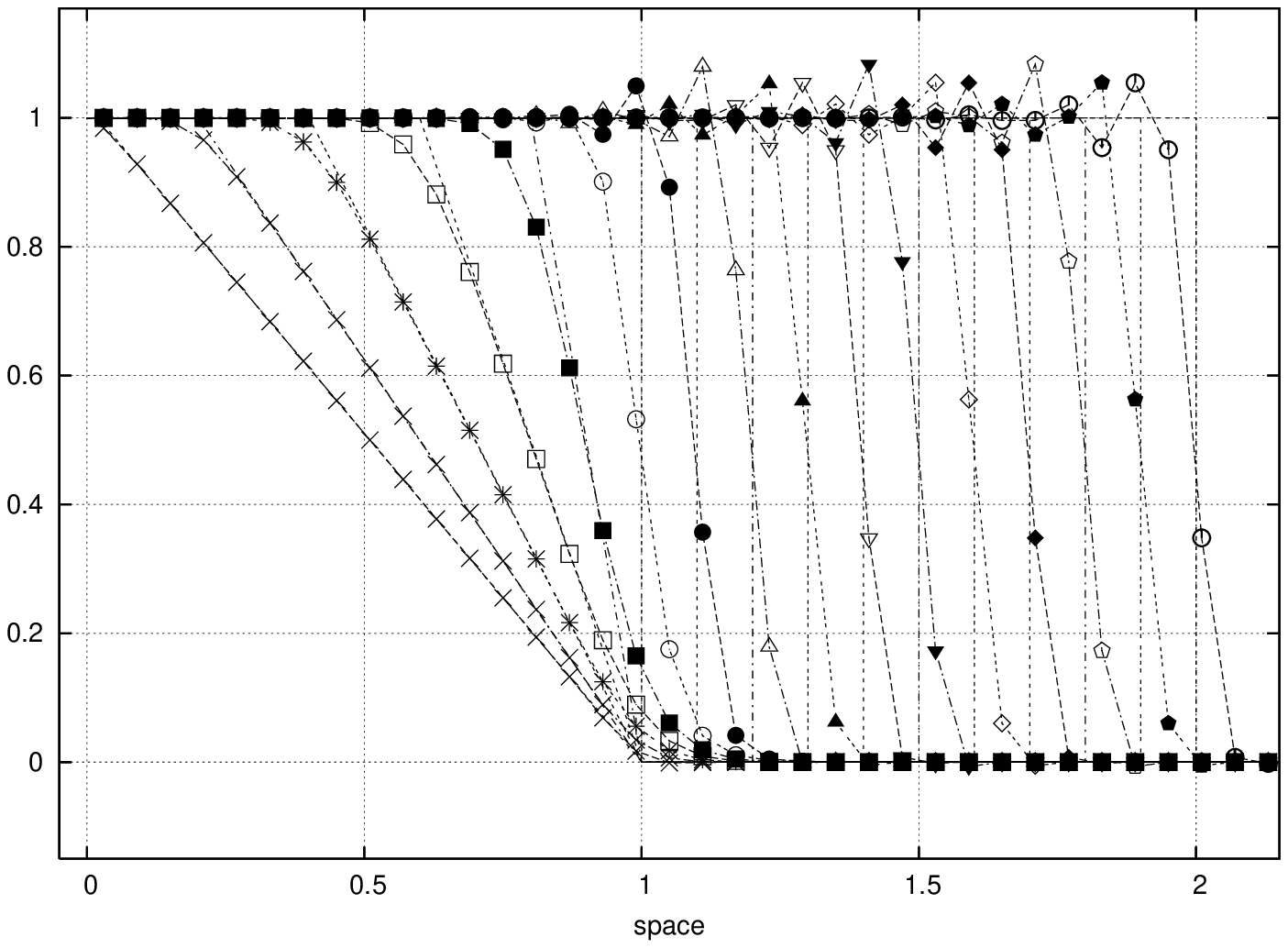} } 
  \bigskip     \noindent  {\bf Figure 6}. \quad 
  Burgers equation. Stable   D1Q3  lattice Boltzmann simulation for 
 a converging shock with equilibrium (\ref{centre-d1q3-f}) 
  associated to the parameters $ \alpha = {1\over2} $,   $  \lambda = 3 $ and 
  $ s_2 = s_3 = 1.7 $.   Computed values are displayed  every 10 time steps.  
   \smallskip  

 \smallskip     \bigskip     \bigskip   
\centerline { \includegraphics[width=.50    \textwidth, angle=-90   ] 
{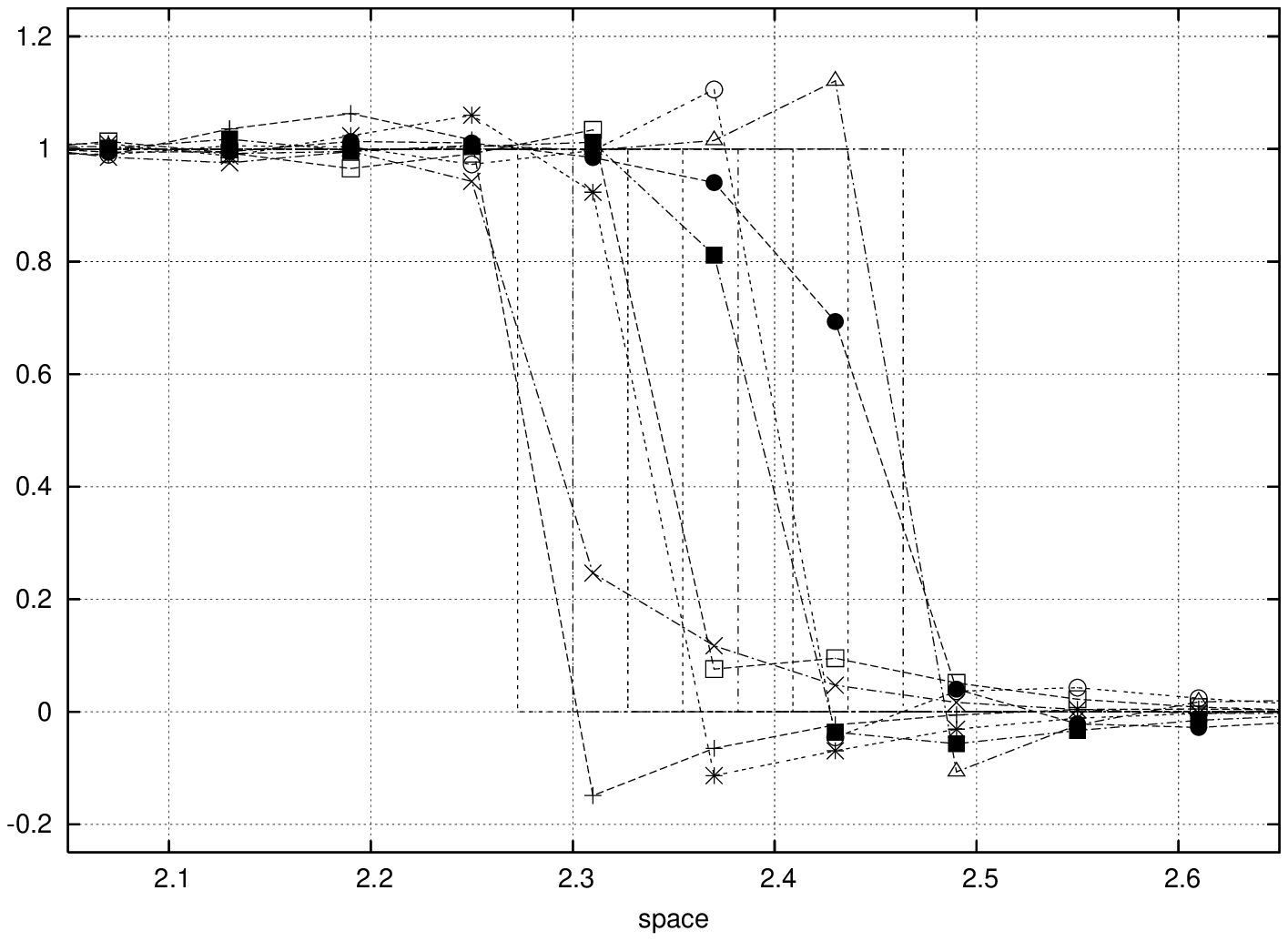}} 
  \bigskip    \noindent  {\bf Figure 7}. \quad  
Burgers equation. Stable   D1Q3  lattice Boltzmann simulation for 
 a converging shock with upwind equilibrium (\ref{upwind-f})
with   $\lambda = 1.1$ and    $ s_2 = s_3 = 1.7 $.  Eight consecutive discrete time steps.  
 \smallskip  

 \bigskip     \bigskip    
\monitem   With the same initial condition (\ref{profil-decroi}), 
we use  the   D1Q3  upwind  version    (\ref{upwind-f}) of the 
lattice Boltzmann scheme. Now the stability condition is not as  severe as in the previous
case and we take   $ \, \lambda = 1.1 $. 
The results, presented in Figure~7, 
are qualitatively analogous 
to the previous one (see Figure~6). 
We observe on  Figure~7 an alternance of monotonic and over or undershooting 
discrete shock profiles.  
%

\monitem    With the same decreasing initial  condition (\ref{profil-decroi}), 
using  the   D1Q2   version    (\ref{distri-equil-d1q2}) leads  
to results   presented on Figure~8. We 
observe only an over-shooting at the discrete shock profile without any 
under-overshooting. 

\smallskip   \bigskip     \bigskip   
\centerline { \includegraphics[width=.50    \textwidth, angle=-90  ]
{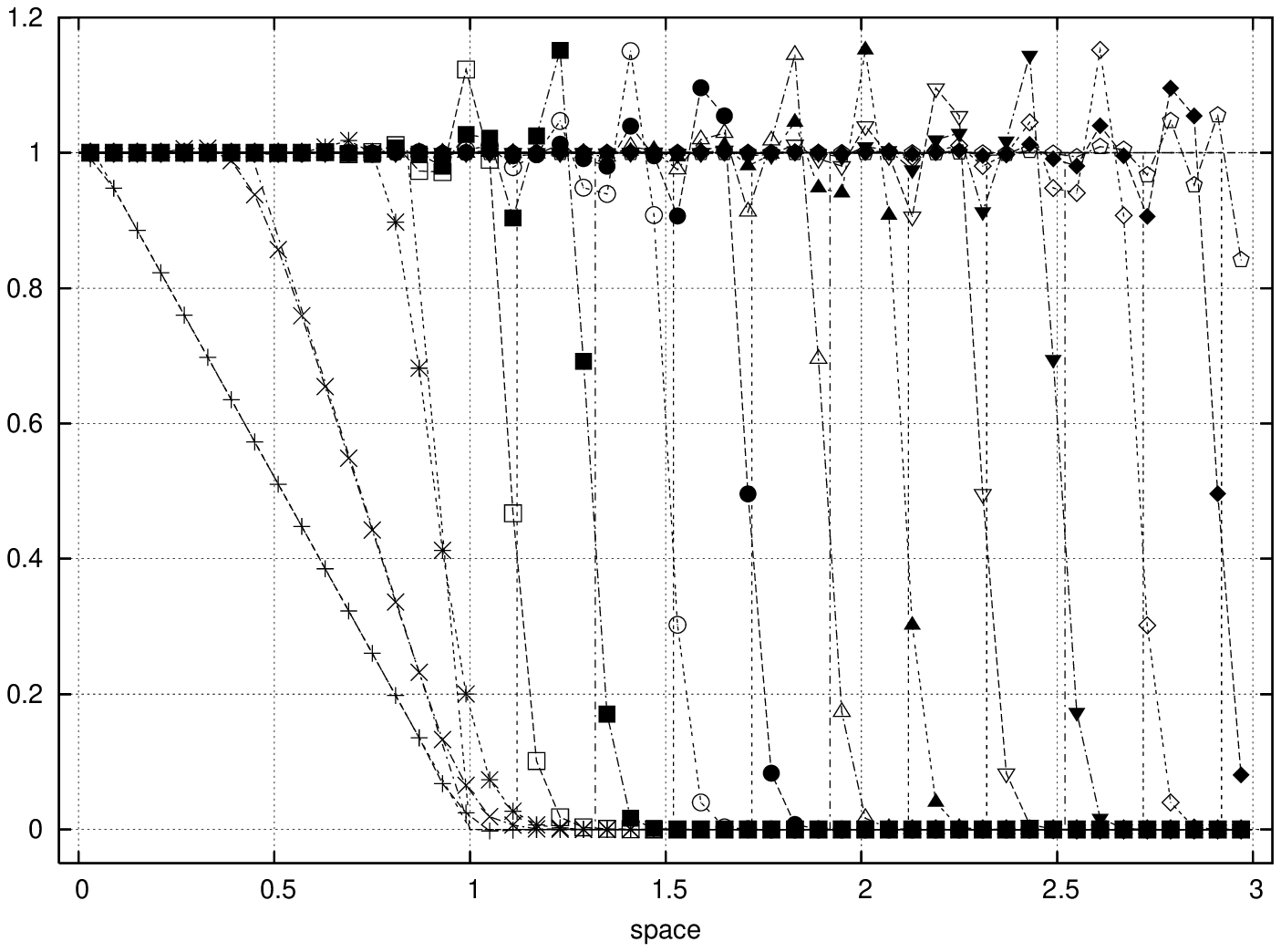}}
  \bigskip  \noindent  {\bf Figure 8}. \quad 
Burgers equation. Stable   D1Q2  lattice Boltzmann simulation for 
 a converging shock with  equilibrium (\ref{distri-equil-d1q2}), 
   $ \lambda = 1.5 \, $  and    $ s_2 = 1.7 $.    
  Computed values are displayed  every 10 time steps. 
\smallskip  

~

\bigskip   \monitem   In a second set of experiments, we  use  the very simple
``two steps'' or ``Riemann''  initial condition. The first one is simply 
\moneq \label{profil-inchoc}    u_0(x)  = 
\left\{  \begin{array} [c] {cl}  \displaystyle  1 \quad & {\rm if}  \qquad   x \,< \, 0.2 \\
0  \quad  & {\rm if} \qquad   x \, > \, 0.2  \, .  \end{array} \right.   
\monend 
The entropic solution of this Riemann problem composed by the 
Burgers equation (\ref{abs-01}) associated with the initial condition 
(\ref{profil-inchoc}) is a discontinuity propagating at the velocity 
$ \, \sigma = {1\over2} \, $   (see {\it e.g.} 
 \cite{GR96},   \cite{DD05} or \cite{La06}). 
With the numerical  schemes introduced previously, this entropy satisfying solution 
is captured with a precision comparable to finite-volume type methods
%
except that for a moving shock, 
a total variation diminishing scheme would not show oscillations ahead and
behind the shock. 
%
The results are presented  on Figure~9  for  numerical schemes
(\ref{centre-d1q3-f}), (\ref{upwind-f}) and  (\ref{distri-equil-d1q2}).
On Figure~10, a zoom of the previous data shows that this
moving  shock is captured by a stencil of four to five mesh points.


\smallskip   \bigskip     \bigskip      
\centerline { \includegraphics[width=.50    \textwidth, angle=-90 ]
{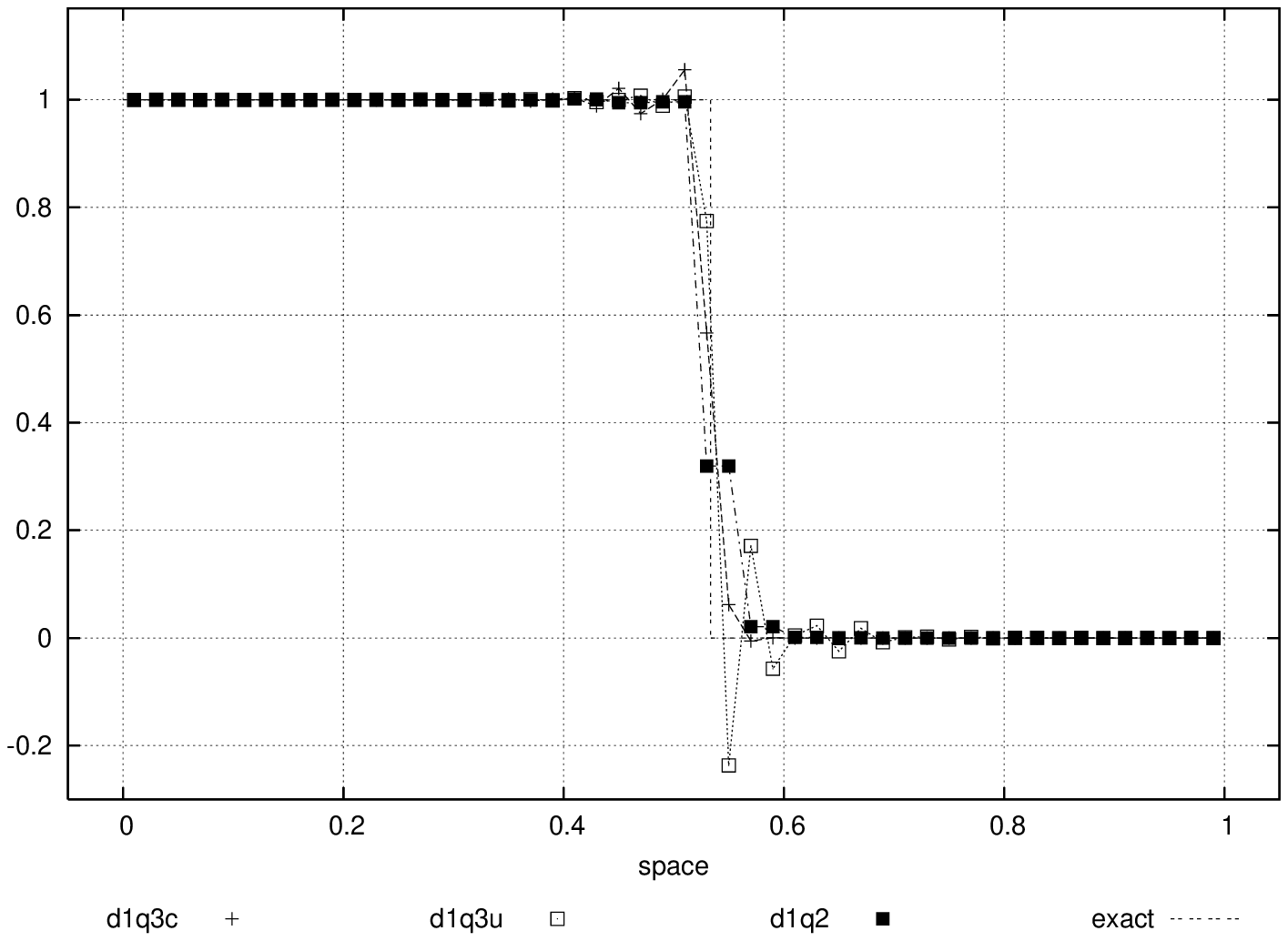}}
\bigskip  \noindent  {\bf Figure 9}. \quad 
The Riemann  problem for the Burgers equation 
associated with the initial condition (\ref{profil-inchoc}) 
develops a shock wave. The figures shows the numerical solutions 
with the three variants of the scheme after 100 discrete time steps
and   parameters  $ \, \lambda = 3  \, $ and $ \, s_2 = s_3 = 1.7 . \,$   
\smallskip  

~
 
\smallskip   \bigskip     \bigskip       
\centerline { \includegraphics[width=.50    \textwidth, angle=-90 ]
{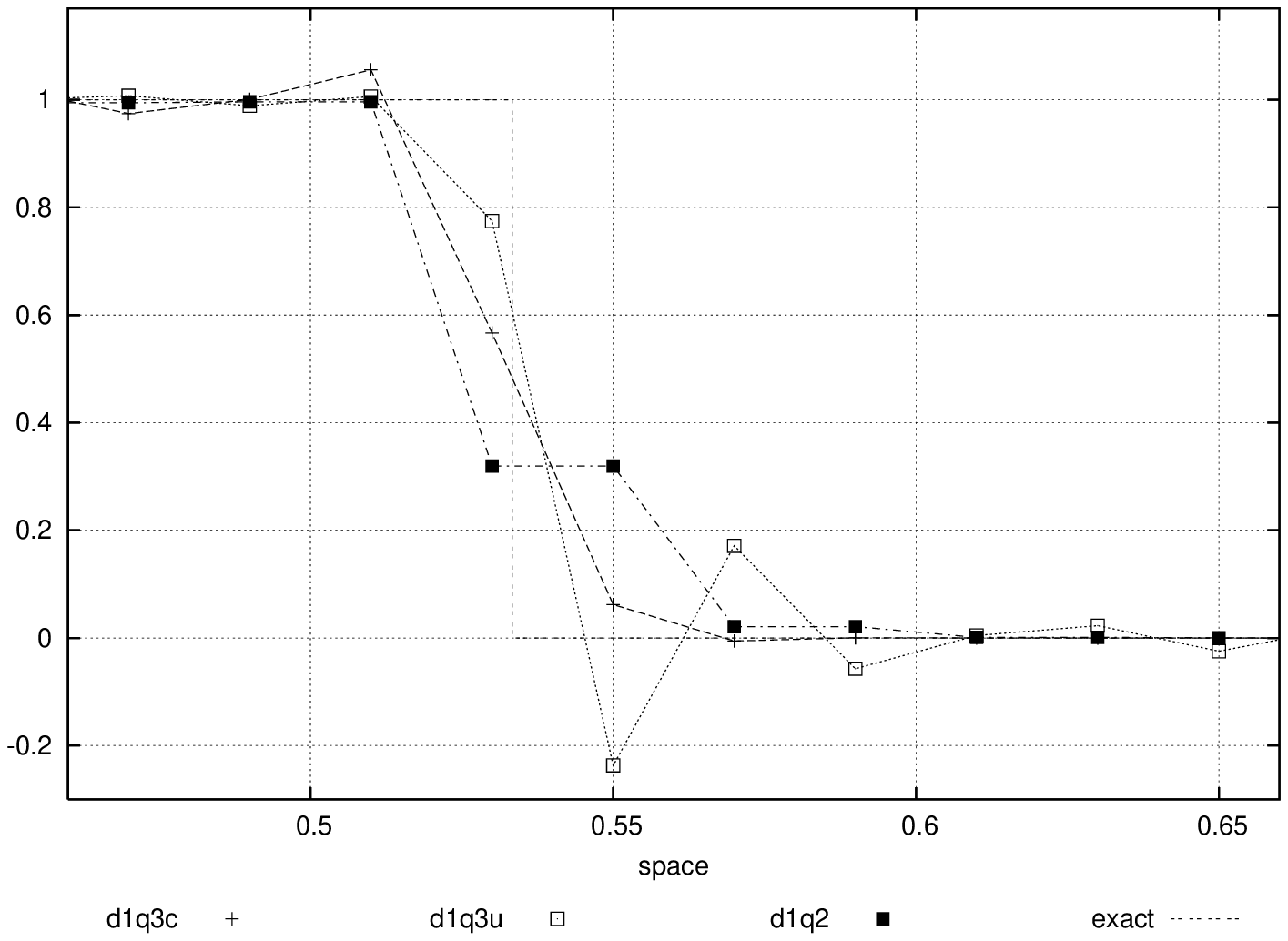}}
\bigskip    \noindent  {\bf Figure 10}. \quad 
Zoom of Figure 10 around the location of the shock wave. 
\smallskip   

\smallskip   \bigskip     \bigskip        
\centerline { \includegraphics[width=.50    \textwidth, angle=-90  ]
{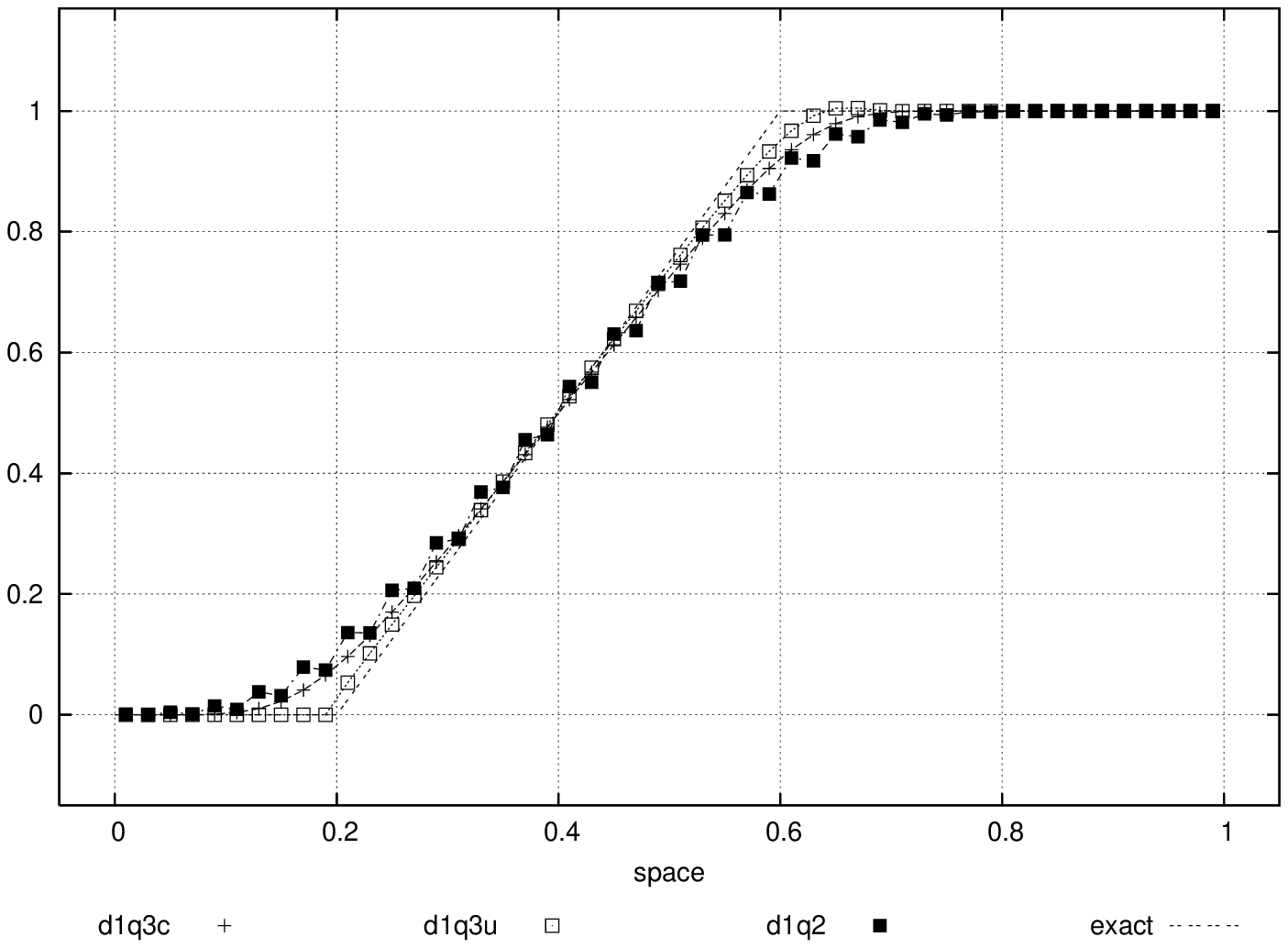}}
\bigskip  \noindent  {\bf Figure 11}. \quad 
The Riemann  problem for the Burgers equation 
associated with the initial condition (\ref{profil-indetente}) 
develops a rarefaction wave. Numerical solutions 
with the three variants of the lattice Boltzmann scheme 
after 100 discrete time steps
and  parameters  $ \, \lambda = 3  , \, $ $ \, s_2 = s_3 = 1.7 . \, $    
\smallskip \smallskip  

\monitem  We reverse the values 0 and 1 in the initial condition (\ref{profil-inchoc})
and obtain in this way a new initial condition~:  
\moneq \label{profil-indetente}    u_0(x)  = 
\left\{  \begin{array} [c] {cl}  \displaystyle  0 \quad & {\rm if}  \qquad   x \,< \, 0.2 \\
1  \quad  & {\rm if} \qquad   x \, > \, 0.2  \, .  \end{array} \right.  
\monend 
The entropic solution  of  (\ref{abs-01})(\ref{profil-indetente}) is
a rarefaction wave~: a continuous solution with two constant states and 
a self-similar component  as  detailed  {\it e.g.}  \cite{GR96},   \cite{DD05} or
\cite{La06}. 
Without any modification of the scheme, the numerical solution with the three previous
variants are presented on Figure~11. 
At  the tricky zones of the foot (Figure~12)
and the top (Figure~13) of the rarefaction,  the slope is discontinuous 
and the solution of the problem   (\ref{abs-01})(\ref{profil-indetente}) 
 is just continuous. We observe that the ``D1Q2'' version of the lattice Boltzmann scheme
 exhibits a two point discrete structure~; in some sense the little number of mesh points 
of this version (\ref{distri-equil-d1q2})  induces some rigidity in the discrete approximation.

\monitem
In this section relative to test cases for unstationary solutions of
the Burgers equation, we have observed two facts. 
First,  
if the dual entropy approach is achieved, the resulting scheme 
is naturally stable even in circumstance where the 
classic linear analysis is {\it a priori} in defect. 
A precise analysis of the competition between nonlinear equilibrium 
and over-relaxation step (\ref{relax-moments})   
%
can be found  the work of  Brownlee {\it et al.} \cite  {BGL07}
with a totally different point of view.    
%
Second, under the convexity condition of the $\, h_j^* \,$ functions 
of the particle decomposition  (\ref{conds-burgers}), 
we observe that the entropy condition is automatically enforced. 
No  so-called rarefaction shock has never been observed with 
the initial condition (\ref{profil-indetente}).

\smallskip  
\centerline { \includegraphics[width=.50    \textwidth, angle=-90  ]
{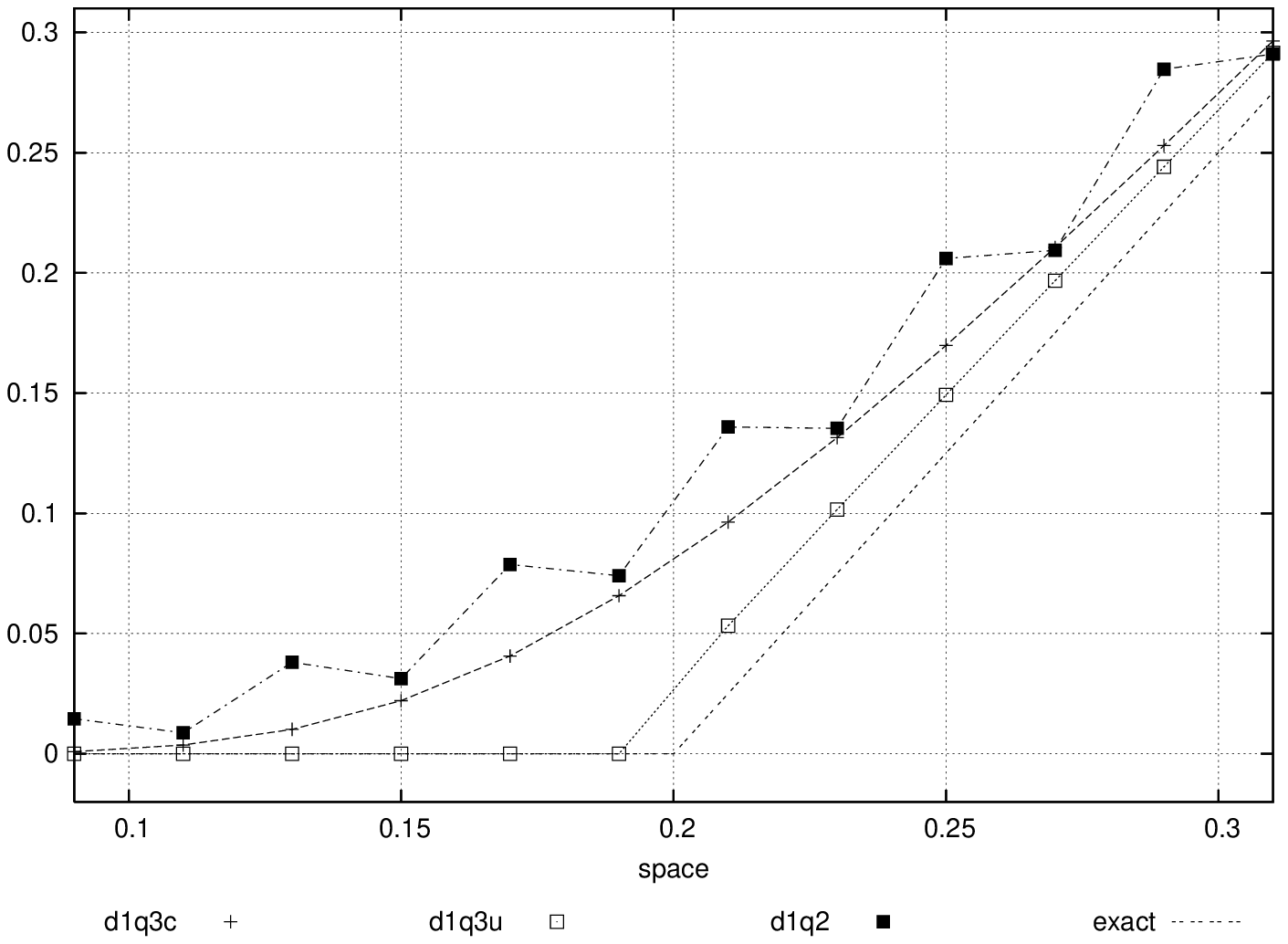}}
\bigskip  \noindent  {\bf Figure 12}. \quad 
Zoom of Figure 12 at the  foot of the rarefaction.   
\bigskip 

\smallskip  
\centerline { \includegraphics[width=.50    \textwidth, angle=-90  ]
{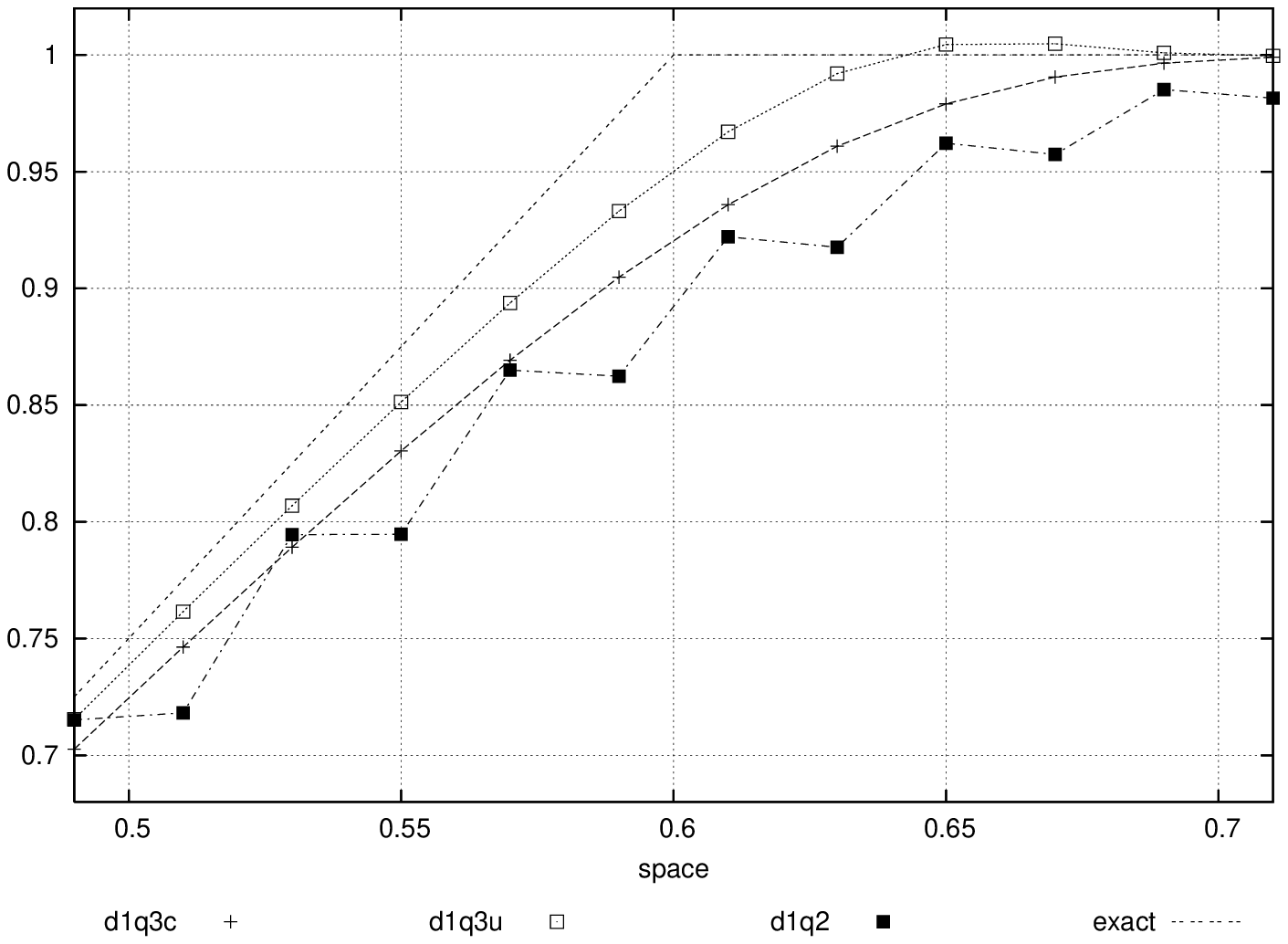}} 
\bigskip     \noindent  {\bf Figure 13}. \quad  
Zoom of Figure 12 at the  top of the rarefaction.  
\smallskip  

\bigskip \bigskip  \noindent {\bf \large 6) \quad 
Linear and nonlinear acoustics }    

\quad 
The extension of the previous ideas from scalar equation to 
hyperbolic  systems is a difficult task. 
We study in this section the first order systems of linear and nonlinear acoustics.

\smallskip \noindent $\bullet$ \quad 
Consider the example of one-dimensional 
linear acoustics with D1Q3 lattice Boltzmann scheme 
to fix the ideas. We recall that we can write 
this physical model as  a hyperbolic system of first order~: 
\moneq  \label{acoust} 
 \partial_t  \begin{pmatrix} \rho \cr q \cr  \end{pmatrix} 
  +  \partial_x \begin{pmatrix} q \cr c_0^2 \, \rho \cr  \end{pmatrix}  = 0   \,.  
\monend 
Then a  mathematical entropy is simply a quadratic form that corresponds to the 
physical energy~: 
\moneq  \label{entrop-acoust} 
 \eta(W) \, \equiv \, 
{{\rho^2}\over{2}}   +  {{q^2}\over{2\, c_0^2}}   \, .
\monend 
The entropy variables are the gradients of the entropy (\ref{entrop-acoust})
relative to the conserved variables $ \, (\rho, \, q) \,$ and we have 
\moneq  \label{entrop-varentrop} 
\varphi = \Big( \rho \,,  {{q}\over{c_0^2}} \Big) \, . 
\monend 
The associated entropy flux  $ \, \zeta(W) \,$ is easy to determine and 
$ \, \zeta(W) =  \rho \, q . \,$ The dual entropy 
$\, \eta^* (\varphi) \equiv \varphi \smb W - \eta(W) \,$ 
and the dual entropy flux  
$\, \zeta^* (\varphi) \equiv \varphi \smb F(W) - \zeta(W) \,$ 
can be evaluated without difficulty and we obtain 
\moneq  \label{entrop-dual} 
\eta^* (\varphi) = \eta(W) \,, \quad  \zeta^* (\varphi) = \zeta(W) \,; 
\monend 
all is quadratic in this system !

\monitem  
We approach the system (\ref{acoust}) with a D1Q3 lattice Boltzmann scheme.  
We use the moments $m$ associated with the same matrix $M$ used 
 for the Burgers equation (see (\ref{abs-04})).  
The associated particle components of the entropy variables 
 $\, \varphi \smb M_j \,$  introduced in (\ref{entroparticle})
 are given according to 
\moneq  \label{acou-entroparticle}  
 \varphi \,\smb \, M_+ \, \equiv \, \rho + {{\lambda \, q}\over{c_0^2}}  \,, \quad 
 \varphi \,\smb \, M_0 \, \equiv \, \rho   \,, \quad 
 \varphi \,\smb \, M_- \,  \equiv \, \rho - {{\lambda \, q}\over{c_0^2}} \, . 
\monend   
The identities (\ref{rep-entrop}) take now the form 
\moneq  \label {eqs-quad} 
\left\{  \begin{array} [c] {cl}  \displaystyle   
 h_+^* \big(  \varphi \,\smb \, M_+ \big) + 
h_0^*  \big(  \varphi \,\smb \, M_0 \big) + 
 h_-^* \big(  \varphi \,\smb \, M_- \big) \, & \equiv\, 
 \displaystyle   \eta^* ( \varphi ) \\ \displaystyle   
 \lambda \, h_+^* \big(  \varphi \,\smb \, M_+ \big) -  \lambda \,
 h_-^* \big(  \varphi \,\smb \, M_- \big) \, & \equiv \, 
 \displaystyle   \zeta^* ( \varphi ) \, .  \end{array} \right. 
\monend   
%
%
We search a  possible solution of system (\ref{eqs-quad})  
 with simple quadratic functions~:
 $ \, h_0^* (y) \equiv a \, y^2 \,$ and 
 $ \, h_+^* (y) \,=\,  h_-^* (y)  \equiv b \, y^2 . \, $
After some lines of algebra, the previous representation and 
the above conditions   (\ref{eqs-quad}) leads  to 
\moneq   \label{acou-hstar}    
\left\{ \begin{array} [c]{cl}  
\displaystyle    h_+^* \Big( \rho + {{\lambda \, q}\over{c_0^2}} \Big)  \, & = \,\, 
\displaystyle  {{ c_0^2}\over{4\, \lambda^2}} \, 
\Big( \rho + {{\lambda \, q}\over{c_0^2}} \Big)^2
\\ \displaystyle   h_0^* (\rho) \, & = \,\, \displaystyle 
 {1\over2} \Big( 1 - {{c_0^2}\over{\lambda^2}}  \Big) \, \rho^2 
\\ \displaystyle  h_-^* \Big( \rho - {{\lambda \, q}\over{c_0^2}} \Big)  \, & = \,\, 
\displaystyle  {{ c_0^2}\over{4\, \lambda^2}} \, 
\Big( \rho - {{\lambda \, q}\over{c_0^2}} \Big)^2 \, . 
\end{array} \right.  
\monend 
The  functions proposed in (\ref{acou-hstar}) are convex under the  stability 
condition~: 
\moneq    \label{acou-stabi}    
 \mid \! c_0 \! \mid \, \leq \, \lambda \, . 
\monend  
This inequality means that 
the numerical waves go faster than the physical ones, a familiar interpretation
of the Courant-Friedrichs-Lewy condition (see {\it e.g.} \cite{La06}). 
A microscopic entropy   $ \, H(f) \,=\, h_+(f_+) + h_0(f_0) + h_-(f_-)  $ 
can be easily derived from (\ref{acou-hstar}) with the following contributors~: 
\moneqstar  
 h_+ \big( f_+ \big)  \,  = \,  {{ \lambda^2}\over{c_0^2}} \,   f_+ ^2  \,, \quad 
h_0 \big( f_0 \big)  \,  = \,  {{1}\over{2 \,  
\Big( \displaystyle  1 - {{c_0^2}\over{\lambda^2}}  \Big)}}  \, f_0 ^2  \,, \quad 
  h_- \big( f_- \big)  \,  = \,  {{ \lambda^2}\over{c_0^2}} \,   f_- ^2   \, . 
\monendstar
The particle distribution 
 $ \,  f_j^{\rm eq} \,$ at equilibrium is a direct consequence of  relations
(\ref{flux-eq}) and (\ref{acou-hstar}) and we have
\moneq   \label{acou-fequil}
 f_+^{\rm eq} \,=\, {{ c_0^2}\over{2\, \lambda^2}} \, 
\Big( \rho + {{\lambda \, q}\over{c_0^2}} \Big) \,, \quad  
 f_0^{\rm eq} \,=\, {1\over2} \Big( 1 - {{c_0^2}\over{\lambda^2}}  \Big)  \, \rho  
\,, \quad  
 f_-^{\rm eq} \,=\, {{ c_0^2}\over{2\, \lambda^2}} \, 
\Big( \rho - {{\lambda \, q}\over{c_0^2}} \Big)  \, . 
\monend 
In terms of moments, the relations (\ref{acou-fequil}) reduce to 
$ \, m_3^{\rm eq} \,=\, c_0^2 \, \rho \,$ as 
%
proposed in  Qian {\it et al.}  \cite{QHL92}.  
Observe that the  equilibrium  (\ref{acou-fequil})  for acoustics
satisfies the dual entropy approach if 
the CFL condition (\ref{acou-stabi}) is satisfied.  


\monitem  
We propose now to introduce a  system of nonlinear acoustics 
  obtained by replacing the 
linear pressure law  in (\ref{acoust}) by a nonlinear one. 
We consider to fix the ideas the particular example of barotropic 
pressure law  $ \, p(\rho) \,$ given according to 
\moneq  \label{nl-pressure} 
p(\rho) = {{1}\over{\gamma}} \, \rho_0 \, c_0^2 \, \Big( {{\rho}\over{\rho_0}} \Big)^\gamma \, ,
\monend
with $ \, \gamma > 1 . \,$  
The corresponding nonlinear system of equations is quite similar to the 
so-called {\it p-system}. It can be written as 
\moneq  \label{nl-acoust} 
 \partial_t  \rho  +  \partial_x  q = 0 \,, \quad 
 \partial_t  q  +   \partial_x  \big( p(\rho) \big) = 0  \,.  
\monend 
It admits a mathematical entropy $\, \eta \,$ and an associated 
entropy flux $\, \zeta \,$ satisfying 
\moneq  \label{nl-acoust-entropy} 
\eta (W)  \, = \, \Phi(\rho) + {{q^2}\over{2}} \,, \quad 
\zeta (W)  \, = \, p(\rho) \, q   \,,  
\monend 
where $ \, \Phi(\smb) \,$ is a primitive of the function 
 $ \, p (\smb) \,$ introduced at the relation (\ref{nl-pressure}). 
In consequence of (\ref{nl-acoust-entropy}), 
the entropy variables $\, \varphi \equiv (\alpha, \, \beta) \,$ take the form 
\moneq  \label{nl-acoust-entropy-variables} 
\alpha \,=\, p(\rho) \,, \quad \beta \,=\, q \, .
\monend 
The dual entropy $\, \eta^* (\smb) \,$ and dual entropy flux  $\, \zeta^* (\smb) \,$
admit the expressions 
\moneq  \label{nl-acoust-dual-entropy} 
\left\{ \begin{array} [c]{cl}   \displaystyle 
\eta^* (\alpha, \, \beta)   \,&=\,  \displaystyle 
 {{\rho_0^2 \, c_0^2}\over{\gamma+1}}  \, \Big( {{\gamma \, \alpha}
\over{\rho_0\, c_0^2}} \Big)^{{\gamma+1}\over{\gamma}} 
 + \,  {{\beta^2}\over{2}} \, \equiv \, 
 {{\rho_0^2 \, c_0^2}\over{\gamma+1}} \, \Big( {{\rho}\over{\rho_0}} 
\Big)^{\gamma+1}  + \, {{\beta^2}\over{2}} \\
\zeta^* (\alpha, \, \beta)   \,&=\, \displaystyle 
  \alpha \, \beta \, \equiv  \, \zeta(\rho,\, q) \, .  
 \end{array} \right.  \monend   

\monitem  
With the matrix $M$ introduced at relation (\ref{abs-04}), 
we denote by $ \, \varphi_+ ,\, $ 
$ \, \varphi_0 \, $ and  $\, \varphi_- \, $ 
the particle components of the entropy variables  
$ \, \varphi \smb M_j \,$ and we have  
\moneq   \label{d1q3-system-particle-components} 
\varphi_+ = \alpha + \lambda \, \beta \,,\quad 
\varphi_0 = \alpha \,,\quad  
\varphi_- = \alpha - \lambda \, \beta \, . 
\monend
It is possible to find nonlinear convex functions 
satisfying   (\ref{eqs-quad}) with the new entropy data 
(\ref{nl-acoust-dual-entropy}). 
By differentiating the relations (\ref{nl-acoust-dual-entropy})
relative to the two entropy variables (\ref{nl-acoust-entropy-variables}), 
 the equilibrium functions $\, f_+^{\rm eq} $, $\, f_0^{\rm eq} \,$ and $\, f_-^{\rm eq} \, $
must satisfy the relations 
\moneq  \label {nl-acoust-fluxeq} 
\left\{ \begin{array} [c]{cl}   \displaystyle 
 f_+^{\rm eq} ( \alpha + \lambda \, \beta ) + 
 f_0^{\rm eq} ( \alpha  ) + 
 f_-^{\rm eq} ( \alpha - \lambda \, \beta ) 
 \,&=\,  \rho  \\  \displaystyle 
\lambda \, f_+^{\rm eq} ( \alpha + \lambda \, \beta ) - 
\lambda \, f_-^{\rm eq} ( \alpha - \lambda \, \beta )  
 \,&=\,  q \equiv \beta  \\  \displaystyle 
\lambda^2 \, f_+^{\rm eq} ( \alpha + \lambda \, \beta ) +
\lambda^2 \, f_-^{\rm eq} ( \alpha - \lambda \, \beta ) 
 \,&=\,  p(\rho) \equiv \alpha \, . 
 \end{array} \right.  \monend   
Then 
%
%
\moneq  \label {nl-acoust-fluxeq-2}  
f_+^{\rm eq} ( \alpha + \lambda \, \beta )  \,=\, 
 {{1}\over{2\, \lambda^2}} \,  ( \alpha + \lambda \, \beta ) \,, \,\, 
 f_0^{\rm eq} ( \alpha )  \,=\, \displaystyle  \rho - {{\alpha}\over{\lambda^2}} \,, \,\, 
 f_-^{\rm eq} ( \alpha - \lambda \, \beta )  \,=\, 
{{1}\over{2\, \lambda^2}} \,  ( \alpha - \lambda \, \beta )  
\monend
and by integration of (\ref{flux-eq}) and  (\ref{nl-acoust-fluxeq-2}), 
we deduce  that the  relations (\ref{acou-hstar}) have to be replaced by 
\moneq  \label{nl-acoust-hstar}
  h_+^* (\alpha) =  h_-^*  (\alpha) = {{1}\over{4\, \lambda^2}} \, \alpha^2 \,, \quad 
  h_0^* (\alpha) =  
 {{\rho_0^2 \, c_0^2}\over{\gamma+1}}  \, \Big( {{\gamma \, \alpha}
\over{\rho_0\, c_0^2}} \Big)^{{\gamma+1}\over{\gamma}}  - \, {{\alpha^2}\over{2 \, \lambda^2}} 
\,.  \monend 
%
The function $ \,   h_+^* (\smb) \equiv   h_-^* (\smb) \, $ is
clearly convex and it is also the case for the function 
 $ \,   h_0^* (\smb) \, $ if its second derivative relative to $ \, \alpha \, $ is
positive, {\it id est} if and only if the following ``dual
stability condition''  is satisfied: 
\moneq  \label{nl-acoust-stabi}
\Big({{\rho}\over{\rho_0}} \Big)^{\gamma-1} \, \Big( {{c_0}\over{\lambda}} \Big)^2  \, \leq \, 1 \, . 
\monend 
%

\smallskip   
\centerline { \includegraphics[width=.50 \textwidth, angle=-90  ]
{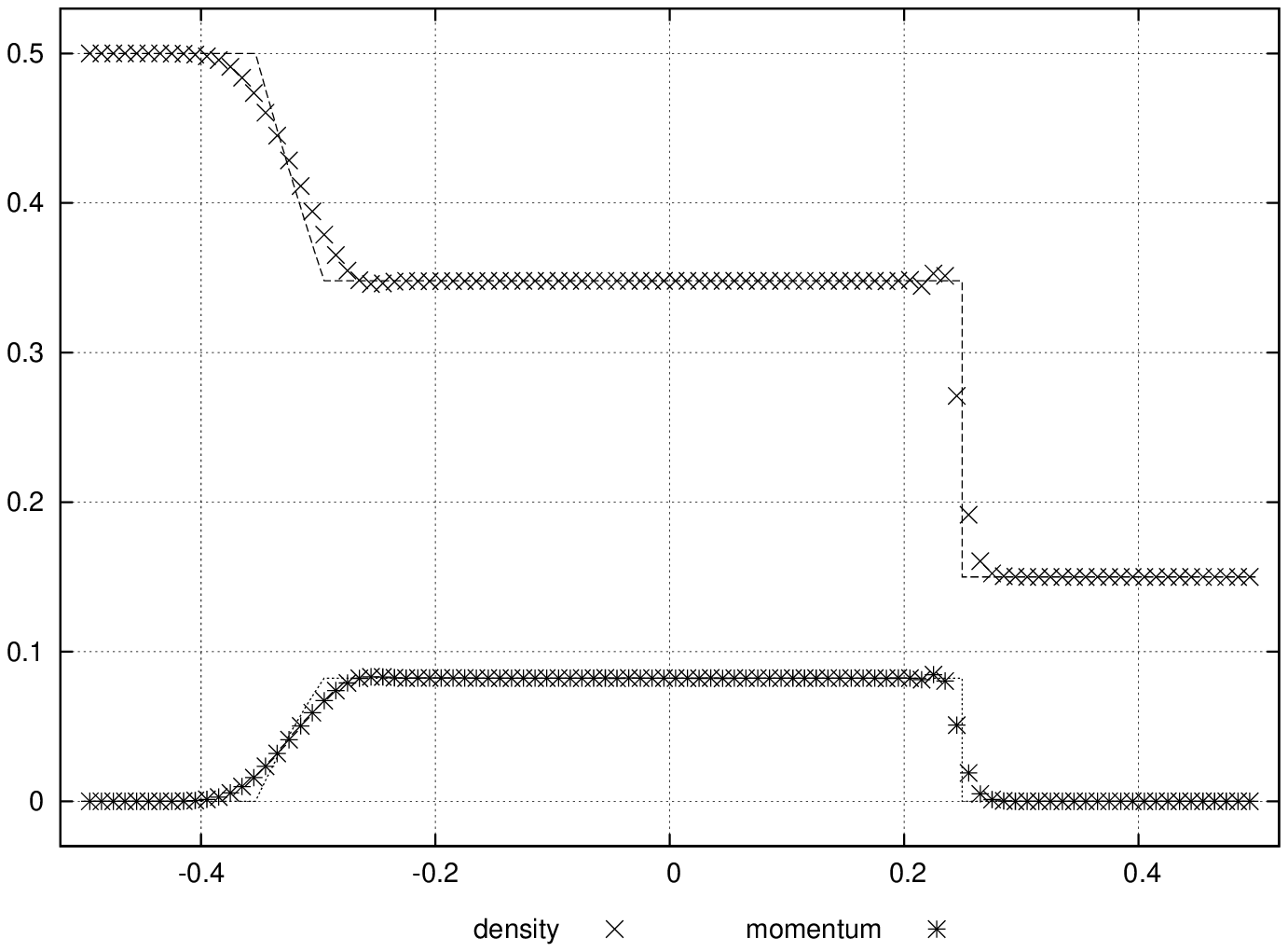} } 
 \smallskip   \smallskip \noindent  {\bf Figure 14}. \quad  
Riemann problem  (\ref {nl-acoust})  (\ref {nl-acoust-riemann-ci})   for the 
system of nonlinear acoustics. The numerical data are precised at the relations
  (\ref{nl-acoust-physical-data}). A rarefaction wave is propagating from
right to left and a shock wave from left to right. Exact (dotted lines) 
and approximated (discrete symbols) profiles of 
density (top) and  momentum (bottom) for 100 mesh points and 60 time steps.   
 \smallskip   \smallskip  

\monitem  
We have tested the system of nonlinear acoustics 
(\ref{nl-pressure}) (\ref{nl-acoust}) with a D1Q3 lattice Boltzmann scheme
for a Riemann problem. 
The initial condition is a discontinuity at $x=0$:
\moneq  \label {nl-acoust-riemann-ci} 
( \rho(x,\, 0) \,,\, q(x, \, 0) ) \,=\, 
\left\{ \begin{array} [c]{ll}   \displaystyle  
( \rho_\ell  \,,\, q_\ell   ) \, \,&{\rm if} \, \, x < 0  \\  \displaystyle 
( \rho_r  \,,\, q_r  ) \, \,&{\rm if} \, \, x > 0 \,. 
 \end{array} \right.  \monend   
We have chosen the physical and numerical parameters
as follows: 
\moneq  \label {nl-acoust-physical-data} 
\gamma =  2  \,,\quad 
{{\rho_\ell}\over{\rho_0}} = 0.5 \,,\quad 
{{\rho_r}\over{\rho_0}} = 0.15 \,,\quad  
q_\ell  =  q_r = 0  \,,\quad  
{{\lambda}\over{c_0}} = 1.2 \,,\quad 
s_3 = 1.7 \, .  
\monend   
The exact solution of the nonlinear hyperbolic system 
 (\ref {nl-acoust})  (\ref {nl-acoust-riemann-ci}) can be obtained without 
difficulty with the general methods presented in     \cite{DD05} or \cite{GR96}.
In the case of  initial data (\ref{nl-acoust-riemann-ci})  (\ref{nl-acoust-physical-data})
a rarefaction wave propagates with a negative velocity and a shock wave
propagates with a positive velocity $ \, \sigma =  0.416 \, \, c_0   . $ 
An intermediate state with 
$\, \rho^* =   0.348  \, \, \rho_0 \, $ and 
$\, q^* =  0.0824  \, \,  \rho_0 \, c_0 \,$ separates these two nonlinear waves. 
With the parameters    (\ref{nl-acoust-physical-data}), 
the condition (\ref{nl-acoust-stabi})  is realized:
$\, \big({{\rho}\over{\rho_0}} \big)^{\gamma-1} \, \big( {{c_0}\over{\lambda}} \big)^2
\leq   0.347 $. 
The numerical results are presented at Figure~14. The rarefaction wave and the shock wave
are correctly captured as in the case of the Burgers equation (see figures 9 and 10). 
When the dual stability condition  (\ref{nl-acoust-stabi}) is not satisfied, 
the lattice Boltzmann scheme replaces the rarefaction by a spurious shock wave 
and becomes completely unusable for higher values of the  parameter defined by the
left hand side of  (\ref{nl-acoust-stabi}). 

\monitem  
As a summary of this section, the generalization of what have been done in this
contribution for the Burgers equation 
with the D1Q3 lattice Boltzmann scheme 
is essentially nontrivial. 
It is possible  to simulate specific nonlinear systems 
of conservation laws and we have experimented with the case
of nonlinear acoustics. 


\bigskip \bigskip  \noindent {\bf \large 7) \quad 
The case of shallow water equations }    

\quad 
The case of shallow water equations has been considered with the lattice Boltzmann scheme 
by Salmon \cite{Sa99} for oceanography applications. 
In the case of   one space dimension  
we can apply the program presented above for   linear and nonlinear acoustic models
and try to represent the dual entropy with the help of a D1Q3 particle distribution. 
We will see in the following the kind of difficulties that we encounter  
with  the dual entropy approach 
with  the present choice  of  a single particle distribution. 

\monitem  
More precisely, we consider the one-dimensional system of conservation laws 
due to  Barr\'e  de Saint Venant~:
\moneq  \label{saintvenant}   
\partial_t  \rho +  \partial_x q \, =  \, 0 \,, \qquad 
\partial_t  q  +   \partial_x \Big( {{q^2}\over{\rho}} + k \, \rho^\gamma \Big)
  \, =  \, 0 \, ,   
\monend 
%
%
where $ k > 0 $ and $ \gamma \geq  1 $ are given positive constants. 
%
We detail in the following the case $ \, \gamma > 1 \,$; 
the case $ \, \gamma \! = \! 1 \,$ is presented in the annex
and conducts to analogous conclusions. 
%
We   introduce  velocity $\, u , \,$   pressure $p$ and sound velocity 
$ \, c>0  \, $ according to  the relations
\moneq  \label{pression-saintvenant}  
 u \equiv {{q}\over{\rho}} \,, \quad p \equiv  k \, \rho^\gamma  \,, \quad
c^2 \equiv   {{\gamma \, p}\over{\rho}} \,=\,   \gamma\, k \, \rho^{\gamma - 1} \, .  
\monend 
Then the entropy $ \, \eta \, $ and the entropy flux $ \, \zeta \,$ satisfy 
\moneq  \label{entrop-stvenant}   
 \eta \,=\, {1\over2} \, \rho \, u^2 + {{p}\over{\gamma - 1}} \,, \quad 
\zeta \,=\, \eta \, u + p \, u \, ; 
\monend  
the entropy variables $\, \varphi = \big(  \partial_\rho \eta ,\, 
 \partial_q  \eta \big) $ $ \, \equiv ( \alpha ,\, \beta ) \,$ can be evaluated without
 difficulty~:
\moneqstar   
\alpha = {{c^2}\over{\gamma - 1}} - {{u^2}\over{2}} \,, \quad \beta = u \, . 
\monendstar
The dual entropy $\, \eta^* \,$ and     the dual entropy flux 
$\, \zeta^* \,$ can be expressed as functions of
the entropy variables~: 
\moneq  \label{entrop-star}  
 \eta^* \,=\, K \, \bigg( \alpha + {{\beta^2}\over{2}} \bigg) 
^{{ \displaystyle \gamma}\over{ \displaystyle \gamma -1}} \,, \quad  
 \zeta^* \,=\, K \, \bigg( \alpha + {{\beta^2}\over{2}} \bigg) 
^{{ \displaystyle \gamma}\over{ \displaystyle \gamma -1}} \, \beta \,, \quad 
K = k \,  \bigg( {{ \displaystyle \gamma -1}\over{ \displaystyle \gamma \, k}}
\bigg)  ^{{ \displaystyle \gamma}\over{ \displaystyle \gamma -1}} \, . 
\monend  
%
We remark that this dual entropy 
$ \, \eta^* \,$ explicited in (\ref{entrop-star})
is no longer the sum of two functions 
of only one entropy variable as in (\ref{entrop-dual}) 
and  (\ref{nl-acoust-dual-entropy}) for linear and nonlinear
acoustics respectively. 
%
The particle components of the entropy variables  
 $ \, \varphi_+ ,\, $ 
$ \, \varphi_0 \, $ and  $\, \varphi_- \, $ 
are still given by the relations (\ref{d1q3-system-particle-components}). 
%
The unknown {\bf convex}  functions $ \, h^*_j \,$ satisfy the identities 
(\ref{eqs-quad}) and  take now the form 
\moneq  \label {eqs-stvenant} 
\left\{ \begin{array} [c]{cl}   \displaystyle 
 h_+^* (  \varphi_+ ) + 
h_0^*   (  \varphi_0 ) +   h_-^*  (  \varphi_- )
 \,&=\,  K \, \bigg(   \displaystyle  \alpha + {{\beta^2}\over{2}} \bigg) 
^{{ \displaystyle \gamma}\over{ \displaystyle \gamma -1}} \\  \displaystyle  
 \lambda \, h_+^*  (  \varphi_+ )   -  \lambda \,
 h_-^*  (  \varphi_- ) \, & = \,  K 
\, \bigg(  \displaystyle   \alpha + {{\beta^2}\over{2}} \bigg) 
^{{ \displaystyle \gamma}\over{ \displaystyle \gamma -1}} \, \beta \, . 
 \end{array} \right.  \monend   

\monitem  
We prove in the following that the system of equations (\ref {eqs-stvenant})
where the unknowns are the convex functions $ \, h_+^* , \,$  $  h_0^* \,$ 
and $ \, h_-^* \,$  of a single real variable, has no solution. 
In order to establish this property, 
we introduce the equilibrium distributions $ \, f_j^{\rm eq} \,$ according to 
({\ref{flux-eq}). We differentiate the relations  (\ref{eqs-stvenant})
relatively to $\, \alpha \,$ and $\, \beta .$ We obtain 
 relations  very similar to (\ref{nl-acoust-fluxeq}):
\moneq  \label {fluxeq-stvenant} 
\left\{ \begin{array} [c]{cl}   \displaystyle 
 f_+^{\rm eq} ( \alpha + \lambda \, \beta ) + 
 f_0^{\rm eq} ( \alpha  ) + 
 f_-^{\rm eq} ( \alpha - \lambda \, \beta ) 
 \,&=\,  \rho  \\  \displaystyle 
\lambda \, f_+^{\rm eq} ( \alpha + \lambda \, \beta ) - 
\lambda \, f_-^{\rm eq} ( \alpha - \lambda \, \beta )  
 \,&=\,  \rho \, u \\  \displaystyle 
\lambda^2 \, f_+^{\rm eq} ( \alpha + \lambda \, \beta ) +
\lambda^2 \, f_-^{\rm eq} ( \alpha - \lambda \, \beta ) 
 \,&=\,  \rho \, u^2 + p \, . 
 \end{array} \right.  \monend   
We are supposed to determine an increasing function $ \,  f_0^{\rm eq} 
\,$ of {\bf only one}  real variable $ \, \alpha \, $  such that 
\moneq  \label {impossible} 
 f_0^{\rm eq} \bigg(  {{c^2}\over{\gamma - 1}} - {{u^2}\over{2}} \bigg) 
 \,\equiv \,  \rho -  {{1}\over{\lambda^2}} \big( \rho \, u^2 + p \big)  \, . 
 \monend 
Due to the elementary calculus 
$ \, {{{\rm d}c^2}\over{{\rm d}\rho}} = \gamma \, k \, (\gamma - 1) \rho^{\gamma -2}
= (\gamma -1) {{c^2}\over{\rho}} $, we differentiate the relation 
(\ref  {impossible}) relative to $ \rho$ 
and independently relatively  to $u$. We obtain 
\moneq  \label {der-impossible} 
 {{c^2}\over{\rho}} \,  \big( f_0^{\rm eq}\big)' (\alpha)  +
 {{1}\over{\lambda^2}} \big( u^2 + c^2 \big) = 1 \,, \quad 
-u \,  \big( f_0^{\rm eq}\big)' (\alpha)  + {{2 \, \rho \, u}\over{\lambda^2}}
= 0 \, . 
 \monend 
We extract the derivative   $ \,  ( f_0^{\rm eq} )' (\alpha) \,$ 
from the second equation of 
(\ref {der-impossible}) and report the result in the first equation. We deduce 
\moneq  \label {final-impossible}  
u^2 + 3\, c^2 = \lambda^2 
 \monend 
and this relation can be correct only for exceptional values of velocity and sound
velocity~! This impossibility is mathematically natural~: it is 
in general not possible to represent
a function of two variables (the right hand side of relation (\ref {impossible}))  
by a simple function of only one variable.

\bigskip \bigskip  \noindent {\bf \large 8) \quad Conclusion and perspectives}   
   
\quad 
We first propose a summary of  the algebraic work that 
%
a ``user'' has to do 
%
 in order to determine  in which domain a given lattice Boltzmann scheme
satisfies the dual stability condition initially  
proposed by Bouchut \cite{Bo03}. 
If very  interesting results are computed with 
 a very good  lattice Boltzmann scheme
 in the framework proposed by d'Humi\`eres \cite{DDH92}, 
%
the procedure follows  five steps. Suppose that the conserved variables 
\moneqstar   
W_k \equiv \sum_j \, M_{kj} \,  f_j  \, 
 \monendstar  
are determined. 
%
Then the convective fluxes follow the relation 
\moneqstar   
  F_{\alpha k}(W)  \equiv \sum_j \, M_{kj} \, v_j^\alpha \,   f_j^{\rm eq} .
 \monendstar   
 First   it is necessary to have a kinetic decomposition 
  of the entropy  and the associated entropy flux of the type 
\moneqstar   
  \eta(W) =  \sum_j  h_j  (f_j^{\rm eq} ) \,, \qquad  
  \zeta_\alpha (W) = \sum_j  v_j^\alpha \,  h_j ( f_j^{\rm eq} ) . 
 \monendstar   
 Second  determine  the   entropy variables  
\moneqstar   
 \varphi = \nabla_W \eta(W) 
 \monendstar   
 and the  one to one mapping  between $W$ and $\varphi$. 
 Third  evaluate  the    Legendre-Fenchel-Moreau duals
\moneqstar   
  h_j^*(y) \equiv  \sup_{f} \, ( y \, f -  h_j (f) )  
 \monendstar   
 of the scalar functions  $  \,  h_j (\smb) . \,  $
  Fourth    determine in which domain    all  the functions 
\moneqstar   
 \varphi \, \longmapsto  h_j^*\big(\varphi\,\smb\, M_j \big)
 \monendstar   
 are     convex. 
    Fifth    report this domain in the $f$ space...

\monitem  
Second, we recall that in this contribution, 
%
we have studied 
 the role of Bouchut stability and convex decomposition of
the dual entropy to develop stable lattice Boltzmann schemes 
in case of simulation of shock and rarefaction waves. 
%
We have applied the above procedure 
 to the Burgers equation,   a fundamental nonlinear scalar equation. 
Then  nonlinear stability does not reduce to a  simple criterion 
 on the relaxation time parameters of the lattice Boltzmann scheme. 
%
A lattice Boltzmann scheme is in general not a finite volume scheme 
and the correct capture of shock waves presented in this contribution
is mathematically absolutely non trivial. 
%
It remains open for us to understand  
why the discrete results with the lattice Boltzmann scheme 
are so well interpreted  in terms of  Bouchut's  theory. 
Moreover, it is a natural question to know 
why the entropy condition is naturally enforced
in the context of nonlinearly stable 
 lattice Boltzmann schemes.

\monitem  
%
Third we have observed that  the   situation for general nonlinear systems is 
not satisfactory. 
Even if all the methodology can be used for a simple nonlinear
system as nonlinear acoustics, it is 
mathematically impossible 
to extend this algebraic  construction to  the familiar nonlinear 
 system of Saint-Venant   equations  one space dimension.  
One idea is to keep the approach as a possible approximation of systems of conservation
laws.   
%
Progress could also result from the use of a {\bf vectorial}  particle distribution 
as initially  proposed by Khobalatte and Perthame in   \cite{KP94} and 
developed by Bouchut  \cite {Bo99} 
for the kinetic finite volume approach. 
Observe that this idea  has been also recognized
as very useful in the lattice Boltzmann community 
for the approximation of thermal fluids
%
and magnetohydrodynamics 
%
as suggested respectively by He, Chen and Doolen  \cite{HCD98}
and Dellar \cite{De02}   and
used by  Peng,    Shu and  Chew \cite{PSC04} among others.

\bigskip \bigskip  \bigskip  \newpage \noindent {\bf \large Acknowledgments}   

\noindent  
The author  thanks  Fran\c cois Bouchut, Benjamin Graille 
and Pierre Lallemand for helpful
discussions during the elaboration of this  work. 
Many thanks also to the 
Institute for Computer Based Modeling in Civil Engineering
of 
Technische Universit\"at Braunschweig (Germany) 
and to the ``LaBS project'' (Lattice Boltzmann Solver,  www.labs-project.org), 
funded by the French FUI8 research program, for supporting this contribution.
Last but not least, 
the author thanks the referees for  very interesting  remarks.
Some of them  have been directly 
incorporated into  the present edition of the article.

\bigskip \bigskip  \bigskip  
\noindent {\bf \large Annex. On shallow water equations 
with $\, \gamma = 1$. }    

\noindent  If $\, \gamma = 1 ,$ we introduce a reference velocity $\, c_* \,$ and replace
the pressure law in (\ref{pression-saintvenant}) by $ \, p =  c_*^2 \, \rho .$ 
Then we introduce a reference density $ \, \rho_* \,$ to express 
in a physically consistent manner 
the algebraic expression a mathematical entropy: 
\moneqstar  
\eta  \,=\,  {{q^2}\over{2 \, \rho}} 
+ c_*^2 \, \rho \, \log {{\rho}\over{\rho_*}} \, . 
\monendstar  
Then  
\moneqstar 
\alpha =  {{\partial \eta}\over{\partial \rho}}  
 = c_*^2 \, \Big( 1 +  \log {{\rho}\over{\rho_*}} \Big) - {{u^2}\over{2}} \,, \quad  
 \beta  =  {{\partial \eta}\over{\partial q}} = u \, . 
\monendstar  
The entropy flux $\, \zeta \, $ is still obtained according to the relation 
(\ref{entrop-stvenant}): $\, \zeta =  \eta \, u + p \, u . \,  $
After some  lines of algebra, the dual entropy $\, \eta^* \equiv 
\alpha \, \rho + \beta \, q - \eta \,$ is equal to 
\moneqstar 
\eta^* =   c_*^2 \, \rho = p = \rho_* \, c_* ^2 \, \exp \Big( {{\alpha + \beta^2 / 2}
\over{c_*^2}} -1 \Big) \, 
\monendstar  
and the dual flux  $\, \zeta^* \equiv 
\alpha \, q + \beta \, ( \rho \, u^2 + p )  - \zeta \,$ is equal to 
$ \, \eta^* \, \beta \,$ as in the case $ \, \gamma > 1 .$ 
Then the relations (\ref {eqs-stvenant}) are generalized without difficulty 
and the identity   (\ref {impossible}) can be now written 
\moneqstar 
 f_0^{\rm eq} \bigg(  c_*^2 \, \big( 1 +  \log {{\rho}\over{\rho_*}} \big) 
 - {{u^2}\over{2}}  \bigg) 
 \,\equiv \,  \rho -  {{1}\over{\lambda^2}} \big( \rho \, u^2 + p \big)  \, . 
\monendstar  
By derivation relative to density and velocity, we get respectively 
\moneqstar  
 {{c_*^2}\over{\rho}} \,  \big( f_0^{\rm eq}\big)' (\alpha)  +
 {{1}\over{\lambda^2}} \big( u^2 + c_*^2 \big) = 1 \,, \quad 
-u \,  \big( f_0^{\rm eq}\big)' (\alpha)  + {{2 \, \rho \, u}\over{\lambda^2}}
= 0 \, .  
\monendstar  
We deduce a necessary relation $ \, u^2 + 3 \,  c_*^2 = \lambda^2  ,\, $  very
close to (\ref {final-impossible}). This relation is satisfied
only for exceptional values of velocity as in the case $ \, \gamma \!>\! 1 .$

\bigskip \bigskip \bigskip  

\medskip


\noindent {\bf \large  References } 

 \vspace{-.4cm} 
 

\begin{thebibliography}{99}
 
\bibitem{ACCD92}
F. Alexander, H. Chen, S. Chen, G.  Doolen. 
``A lattice Boltzmann model for compressible fluids'', 
 {\it Physical Review A},  vol.~{\bf 46},   p~1967-1970, 1992.

\bibitem{BGK}
P.  Bhatnagar, E.  Gross and M. Krook. ``A Model for Collision Processes
  in Gases. I. Small Amplitude Processes in Charged and Neutral One-Component Systems'', 
  {\it    Physical Review}, vol.~{\bf 94}, p.~511-525, 1954.


\bibitem{BL87} 
B. Boghosian, C. Levermore. ``A Cellular Automaton for Burgers's Equation'', 
  {\it  Complex Systems},  vol.~{\bf 1}, p~17-30, 1987.

\bibitem{BLY04} 
B. Boghosian, P. Love, J. Yepez.  ``Entropic Lattice Boltzmann Model for 
Burgers' Equation'',   {\it  Philosophical Transactions of the Royal 
Society A},  vol.~{\bf 362},  p~1691-1702, 2004.
 



\bibitem {Bo99} 
F. Bouchut.  ``Construction of BGK models with a family of kinetic 
entropies for a given system of conservation laws'', {\it Journal of  
Statistical Physics}, vol.~{\bf 95}, p.~113-170, 1999.  



\bibitem {Bo03} 
F. Bouchut.  ``Entropy satisfying flux vector splittings and kinetic BGK models'', 
{\it Numerische  Mathematik},  vol.~{\bf  94},   p.~623-672, 2003.


\bibitem {Bo04} 
F. Bouchut.  ``A reduced stability condition for nonlinear relaxation to conservation laws'', 
{\it Journal of Hyperbolic Differential Equations},  vol.~{\bf  1},   p.~149-170, 2004.

\bibitem {Bo10} 
F. Bouchut. Personal communication, Paris, June 2010. 

 
\bibitem {Br64} 
J. Broadwell. ``Study of a rarefied shear flow by the discrete velocity
method'', {\it Journal of Fluid Mechanics},  vol.~{\bf 19}, 
p.~401-14, 1964. 


\bibitem {BGL07} 
R.  Brownlee, A.  Gorban,  J. Levesleye. ``Stability and stabilization
 of the lattice Boltzmann method'', 
{\it Physical Review E}, vol.~{\bf 75}, 036711, 2007.


\bibitem{BBGG2k} 
J.M. Buick, C.L. Buckley, C.A. Greated, J. Gilbert. 
``Lattice Boltzmann BGK simulation of non-linear sound waves~: the development of a shock
front'',   {\it   Journal of Physics A: Mathematical and General}, 
 vol.~{\bf 33}, p.~3917-3928, 2000. 

\bibitem{Ca80}  
H. Cabannes. {\it The Discrete Boltzmann Equation (Theory and Applications)}, 
Lecture Notes,  University of California at Berkeley, 1980.
Revised edition by H.~Cabannes, R.~Gatignol and L.S.~Luo, 2003 
 available at 
http://henri.cabannes.free.fr/Cours${\rm \underbar{~}}$de${\rm \underbar{~}}$Berkeley.pdf.

 
\bibitem{CLL94}  
G.Q. Chen, C.D. Levermore,  T.P. Liu.
``Hyperbolic conservation laws with stiff relaxation and entropy'', 
 {\it Communications on Pure and Applied Mathematics},  vol.~{\bf 47}, 
p.~787-830, 1994. 


\bibitem{CK09} 
S.  Chikatamarla, I. Karlin, ``Lattices for the lattice Boltzmann method'', 
{\it Physical Review E}, vol.~{\bf 79}, 046701, 2009.



 
\bibitem{De02}  
P. Dellar. 
``Lattice Kinetic Schemes for Magnetohydrodynamics'', 
{\it Journal of Computational Physics},   vol.~{\bf 179},  p.~95-126, 2002.
 
\bibitem{DD05}  
 B. Despr\'es, F. Dubois. {\it  Syst\`emes hyperboliques de lois de conservation~;
 Application \`a  la dynamique des gaz}, 
  Editions de l'Ecole Polytechnique, Palaiseau, 2005. 


\bibitem{DL07}  
Y. Duan,    R.X. Liu. ``Lattice Boltzmann model for 
two-dimensional unsteady Burgers' equation'', 
{\it Journal of Computational and Applied Mathematics},  vol.~{\bf 206}, 
p.~432-439, 2007.   

\bibitem{Du08}       
F. Dubois. ``On lattice Boltzmann scheme, finite volumes and boundary conditions'', 
 {\it Progress in Computational Fluid Dynamics}, vol.~{\bf 8}, p.~11-24, 2008.

\bibitem{El91}     
B. Elton. ``A lattice Boltzmann method for a two-dimensional viscous Burgers equation: 
computational results''.  In {\it Supercomputing'91: Proceedings of the 1991
ACM/IEEE conference on Supercomputing}, p.~242-252. New York, NY, USA: ACM Press, 1991.

\bibitem{ELR93}    
B. Elton, C. Levermore, G.  Rodrigue. ``Convergence of Convective-Diffusive 
Lattice Boltzmann Methods'', 
{\it SIAM  Journal on  Numerical Analysis},  vol.~{\bf 32}, p.~1327-1354, 1995.
 
\bibitem{EO80} 
 B. Engquist, S. Osher. ``Stable and entropy satisfying approximations for 
transonic flow calculations'', {\it Mathematics of Computation}, 
vol.~{\bf  34}, p.~45-75, 1980. 

\bibitem{FL71}    
K.O. Friedrichs, P.D. Lax. ``Systems of Conservation Equations with a Convex 
Extension'', {\it Proc. Nat. Sciences USA}, vol.~{\bf 68}, p.~1686-1688, 1971. 

\bibitem{Ga70}   
R. Gatignol. ``Th\'eorie cin\'etique d'un gaz \`a r\'epartition discr\`ete 
de vitesses'', {\it  Zeitschrift f\"ur Flugwissenschaften},  vol.~{\bf 18}, 
 p.~93-97, 1970.  

\bibitem{GR96} 
E.  Godlewski, P.A.  Raviart. {\it 
Numerical Approximation of Hyperbolic Systems of Conservation Laws}, 
    Springer-Verlag,  New York Inc.,  1996.


\bibitem{Go61}
S. K. Godunov.  ``An interesting class of quasilinear systems'',  
{\it Dokl. Akad. Nauk SSSR},   vol.~{\bf 139},  p.~521-523~, 1961; 
see also  {\it Soviet. Math.},  vol.~{\bf 2}, p. 947-949,  1961.   


\bibitem{HCD98}
X.  He, S. Chen, G.D.  Doolen. ``A Novel Thermal Model 
for the Lattice Boltzmann Method in Incompressible Limit'', 
 {\it Journal of Computational Physics},  vol.~{\bf   146}, 
 p.~282-300, 1998.
 



\bibitem{He87}
M. H\'enon. ``Viscosity of a Lattice Gas'', {\it  Complex Systems},   
vol.~{\bf 1},  p.~763-789, 1987.

\bibitem{DDH92}    
D. d'Humi\`eres. ``Generalized Lattice-Boltzmann Equations'', 
in {\it Rarefied Gas Dynamics: Theory
and Simulations}, vol.~{\bf 159} of {\it AIAA Progress in
Astronautics and Astronautics}, p.~450-458, 1992.  

\bibitem{JKL05}
M. Junk, A. Klar, L.S. Luo. ``Asymptotic analysis of the lattice Boltzmann equation'',
 {\it  Journal of Computational Physics},vol.~{\bf 210},  p.~676-704, 2005. 


\bibitem{KA10}
I.V. Karlin, P. Asinari. ``Factorization symmetry in the lattice Boltzmann method'',
{\it Physica A},   vol.~{\bf  389},  p.~1530-1548, 2010. 


\bibitem{KGSB98}
I.V.  Karlin, A.N. Gorban,  S. Succi and V. Boffi. 
``Maximum Entropy Principle for Lattice Kinetic Equations'',  
{\it  Physical Review Letters},  vol.~{\bf 81}, p.~6-9, 1998. 

\bibitem{KP94}
B.  Khobalatte, B. Perthame. ``Maximum principle on the entropy 
and second-order kinetic schemes'', 
{\it Mathematics of Computation},   vol.~{\bf 62}, p.~ 119-131, 1994. 
 



\bibitem{LL00}  
P. Lallemand, L.-S. Luo. ``Theory of the lattice Boltzmann method: 
Dispersion, dissipation, isotropy, Galilean invariance, and stability'',
{\it Physical Review E}, vol.~{\bf 61}, p.~6546-6562, June 2000.  

\bibitem{La06}  
P.D. Lax. {\it Hyperbolic Partial Differential Equations}, 
American Mathematical Society, 
Courant Lecture Notes, vol.~{\bf  14}, 2006.     


\bibitem{NSC08}  
X. Nie, X. Shan, H.  Chen. ``Thermal lattice Boltzmann model 
for gases with internal degrees of freedom'', 
 {\it Physical  Review E},  vol.~{\bf  77}, p.~035701(R), 2008.

\bibitem{PSC04}  	
Y. Peng,  C. Shu, Y.T. Chew. ``A 3D incompressible thermal lattice 
Boltzmann model and its application to simulate natural convection in a cubic cavity'', 
{\it   	Journal of Computational Physics},  
 vol.~{\bf  193}, p.~260 - 274, 2004. 

\bibitem{Pe90}  
B. Perthame. ``Boltzmann Type Schemes for Gas Dynamics and the Entropy Property'', 
{\it SIAM Journal on Numerical Analysis},  vol.~{\bf 27},  p.~1405-1421,  
1990.
 

\bibitem{PHSS90}  
P.C. Philippi, L.A. Hegele, R. Surmas, D.N. Siebert. 
``From the Boltzmann to the Lattice-Boltzmann equation: beyond BGK collision models'', 
{\it International Journal of Modern Physics C, 
Computational Physics and Physical Computation},  vol.~{\bf 18},  p.~556-565, 2007.  


\bibitem{QHL92} 
 Y. Qian, D. d'Humi\`eres,  P. Lallemand, 
``Lattice BGK Models for Navier-Stokes Equation'', 
{\it Europhysics Letters},  vol.~{\bf 17},  p.~479-484, 1992.


\bibitem{QZ98} 
Y. Qian, Y.  Zhou. ``Complete Galilean-invariant lattice BGK models 
for the Navier-Stokes equation'', {\it  Europhysics  Letters}, 
 vol.~{\bf 42},  p.~359-364, 1998. 

\bibitem{Sa99} 
R. Salmon. ``The lattice Boltzmann method as a basis for ocean circulation modeling'',
 {\it Journal of Marine Research},  vol.~{\bf 57},  p.~503-535, 1999. 




\bibitem{Sh48}   
C. Shannon. ``A mathematical theory of communication'', 
{\it The Bell Labs Systems Technical Journal}, vol.~{\bf 27}, 
p.~379-423  and 623-656, 1948.    

\bibitem{Ye02}  
J. Yepez. 
``Quantum Lattice-Gas Model for the Burgers Equation'', {\it Journal
of Statistical Physics}, vol.~{\bf 107},  p.~203-224, 2002.



\end{thebibliography}
\end{document}